%% file: root.tex
\documentclass[journal]{IEEEtran}
\usepackage{multirow}
\usepackage{algorithm,algpseudocode}
\usepackage{comment}
\usepackage{cancel}

\usepackage{multicol}
\usepackage{amssymb}
\usepackage{mathrsfs}
\usepackage{amsmath}
\usepackage[english]{babel}
\usepackage{float}
\usepackage{accents}
\usepackage{bm}
\usepackage{mathtools}
\usepackage{verbatim}
\usepackage{easy-todo}
\usepackage[standard]{ntheorem}
\usepackage{soul}
\usepackage[normalem]{ulem}

\newcommand{\norm}[1]{\left\lVert#1\right\rVert}

\newcommand{\q}{\mathbf{q}}
\newcommand{\uc}{\mathbf{u}}

\newcommand{\bs}[1]{\boldsymbol{#1}}
\newcommand{\bb}[1]{\mathbf{#1}}

\begin{document}
\title{Robust optimal density control of robotic swarms}
\author{Carlo Sinigaglia, Andrea Manzoni, Francesco Braghin, and Spring Berman 
\thanks{Submitted to IEEE Transactions on Automatic Control on \today. This work was supported by the Italian Ministry of Education, University and Research (MIUR).}
\thanks{C. Sinigaglia and F. Braghin are with the Department of Mechanical Engineering, Politecnico di Milano, Milano, 20156 Italy (e-mail: carlo.sinigaglia@polimi.it; francesco.braghin@polimi.it).}
\thanks{A. Manzoni is with the MOX Department of Mathematics, Politecnico di Milano, Milano, 20133 Italy (e-mail: andrea1.manzoni@polimi.it)}
\thanks{S. Berman is with the School for Engineering of Matter, Transport and Energy, Arizona State University, Tempe, AZ 85287 USA (e-mail: spring.berman@asu.edu)}
}

\maketitle

\begin{abstract}
In this paper, we propose a computationally efficient, robust density control strategy for the mean-field model of a robotic swarm. We formulate a static optimal control problem (OCP) that computes a robot velocity field which drives the swarm to a target  equilibrium density, and we prove the stability of the controlled system in the presence of transient perturbations and uncertainties in the initial conditions. The density dynamics are described by a linear elliptic advection-diffusion equation in which the control enters bilinearly into the advection term. The well-posedness of the state problem is ensured by an integral constraint. We prove the existence of optimal controls by embedding the state constraint into the weak formulation of the state dynamics. The resulting control field is space-dependent and does not require any communication between robots or costly density estimation algorithms. Based on the properties of the primal and dual systems, we first propose a method to accommodate the state constraint. Exploiting the properties of the state dynamics and associated controls, we then construct a modified dynamic OCP to speed up the convergence to the target equilibrium density of the associated static problem. We then show that the finite-element discretization of the static and dynamic OCPs inherits the structure and several useful properties of their infinite-dimensional formulations. Finally, we demonstrate the effectiveness of our control approach through numerical simulations of scenarios with obstacles and an external velocity field.
\end{abstract}


\newcommand*{\QEDA}{\null\nobreak\hfill\ensuremath{\square}}
 

\section{Introduction}
\label{section:intro}
\input{Sections/Introduction.tex}

\section{The Optimal Control Problem}
\label{section:ocp}
\input{Sections/analysis}

\section{Numerical Analysis of the OCP}
\label{section:ocp_num}

\input{Sections/numerical_analysis}

\section{Numerical Simulations}
\label{section:sim}
\input{Sections/numerical_simulations}

\section{Conclusions}
\label{section:conclusion}
\input{Sections/conclusion}

\bibliographystyle{IEEEtran}
\bibliography{IEEEabrv,References}

\end{document}

%% file: Sections/Introduction.tex
\IEEEPARstart{L}{arge-scale} collectives 
 of robots, or robotic swarms, are increasingly finding applications in a variety of tasks, 
 such as search-and-rescue missions, 
infrastructure inspection and maintenance, precision agriculture, 
 and many others \cite{Dorigo2021}. This is in part due to the 
significant decrease in the cost of electronic components over the past few decades, which facilitates the fabrication of 
very large numbers 
 of robots. 
Due to size and economic constraints, the computational power of a single swarm member 
is necessarily limited, which restricts the complexity of its control algorithms. 
From a control-theoretic point of view, the challenge is to synthesize controllers 
that can be implemented on swarms of such robots to produce collective behaviors that achieve specified high-level tasks, in a way that accommodates 
the high dimensionality of the system.   

Classical path planning and control algorithms either do not scale well with the number of robots or do not allow the designer to specify complex high-level objectives. Recently, macroscopic descriptions of swarm dynamics in the form of mean-field models \cite{Elamvazhuthi2020} 
have been used to devise robust path planning algorithms for robotic swarms to perform collective tasks such as coverage and mapping \cite{hsieh2008biologically,elamvazhuthi2018pde}. 
Mean-field models provide a general probabilistic framework that can be used to 
design control algorithms for swarms of agents with stochastic behaviors.
In this framework, swarm tasks 
are specified in terms of macroscopic population dynamics that are described by a mean-field model, and this model is used to derive the robot control policies, which guide the microscopic dynamics of individual robots and drive the swarm to collectively reproduce the macroscopic dynamics in expectation.
A consistent way 
of analyzing the performance of such control policies 
when they are implemented on a finite number of robots has been developed in \cite{martinez2021analysis}. 

In the mean-field setting, the robotic swarm
is represented by a probability density, which is independent of the number of robots, that evolves over space and time according to a 
Kolmogorov forward equation. Finite-dimensional mean-field models consist of a linear system of ordinary differential equations (ODEs) describing the dynamics of a swarm that evolves according to a Markov chain over a finite state space, which consists of
a set of tasks or 
discrete spatial locations (e.g., \cite{berman2009optimized}). One type of infinite-dimensional mean-field model is a linear parabolic advection-diffusion partial differential equation (PDE) governing the space-time dynamics of a swarm that follows a deterministic velocity field perturbed by noise, modeled by a Wiener process, over a continuous state space.

Recently developed control strategies for swarms based on mean-field models often rely on estimation of
the local swarm density at each instant \cite{bandyopadhyay2017probabilistic, eren2017velocity}, which is implemented with
decentralized estimation algorithms that are computationally costly and require inter-robot communication. 
The work \cite{eren2017velocity} synthesizes control laws that are inversely proportional to 
the estimated density, which generates
unphysically high velocity fields when the density is small.
On the other hand, Markov chain-based probabilistic algorithms for swarm density control, such as
\cite{accikmecse2012markov}, do not necessarily  
require 
communication between robots 
and provide interesting self-healing properties. However, 
an {\it a priori} discretization of the state space is 
needed to set up the control problem in the 
Markov-chain framework.
The control laws in \cite{accikmecse2012markov} do not depend on  
the estimated swarm density,
but rather 
are implicit functions of the target swarm density and the state-space discretization. 

In this paper, we propose an optimization-based algorithm for computing a feedback control law that steers a swarm of agents with both advective and diffusive motion towards a target equilibrium density. The control law is defined as the advection field and depends only on space, not on the estimated swarm density or the initial conditions of the swarm, and can be preprogrammed on simple agents without inter-agent communication.
Our approach is similar to \cite{accikmecse2012markov}, but in contrast to this work, which uses a finite-dimensional Markov chain model, we consider an infinite-dimensional mean-field formulation where the equilibrium swarm density is 
the solution of a linear elliptic advection-diffusion PDE, which defines the equilibrium condition
of the corresponding time-dependent parabolic problem.
The feedback control law 
is the solution of a static Optimal Control Problem (OCP) whose state dynamics govern the equilibrium swarm density; 
the control field enters bilinearly into the state dynamics. 

The dynamic version of this problem, where the state dynamics are described by an evolution equation in the form of a linear parabolic PDE, has been studied extensively; see, e.g., \cite{roy2018fokker} for both 
theoretical and numerical treatments of the problem where the control field is null at the boundary, 
and \cite{sinigaglia2021density} for a boundary control application. It is known that the dynamic control problem is controllable 
to every sufficiently smooth target distribution \cite{karthik_bil}. When a dynamic OCP is considered, the resulting control field is inherently open-loop and depends on the initial conditions. On the other hand, \cite{zheng2021transporting} 
stabilizes target distributions using a feedback control law that, 
although robust to external transient perturbations, 
is inversely proportional to the local swarm density, which generates unphysically large control fields when this density is close to zero. 
Moreover, the local density must be
estimated by the agents using a communication mechanism that requires additional computational resources.
A distributed algorithm for density estimation to reduce the computational cost is proposed in \cite{zheng2021distributed}. 


In the optimization-based algorithm that we propose, the OCP is designed to compute
an equilibrium swarm density 
that is as close as possible to a 
target density, which may be non-smooth, while balancing the control expenditure. 
The properties 
of the state operator enable us to show that the optimal control field globally stabilizes the equilibrium density, driving the swarm asymptotically to this density from every initial condition.
In addition, our 
approach is computationally efficient since it entails 
the solution of a static OCP at each iteration of a numerical optimization procedure. For cases where 
the initial swarm density is approximately known, we set up a dynamic OCP which makes use of the static solution. The dynamic OCP is formulated in a way that ensures convergence to the static, globally stabilizing control field.

By embedding the state dynamics in a suitable zero-mean subspace and using methods from functional analysis, we prove the existence of static optimal controls where the control advection field is chosen in an appropriate functional space. Then, exploiting the kernel properties of the state dynamics, we prove a stability theorem for the resulting equilibrium density.

From a computational standpoint, we propose a consistent finite-element discretization 
of both the transient and steady-state density dynamics and show that it preserves particular 
properties of the infinite-dimensional problem. 
Based on the kernel properties of the state matrix of the 
algebraic discretization, we develop a numerical algorithm to efficiently compute the reduced gradient. Then, using arguments analogous to those applied to the static OCP, 
we prove the existence of optimal controls for the associated dynamic OCP  and solve this problem
in a similar way, using the static solution as a warm start. We also demonstrate that we are able to ensure convergence of the dynamic solution to the static solution of the OCP for both the state and control variables. This is a robustness property with respect to uncertainties in the initial conditions and transient disturbances. We validate our approach in numerical simulations that include complex scenarios in which the swarm moves through environments with obstacles and an external velocity field.

The paper is organized as follows. In Section~\ref{section:ocp}, the infinite-dimensional formulation of the OCP is presented and analyzed, and then optimality conditions are derived together with a set of useful properties. In Section~\ref{section:ocp_num}, the properties of the finite-element discretization of the OCP are proved and discussed in detail; a solution algorithm which exploits these properties is then proposed. In Section~\ref{section:sim}, three test cases are solved numerically to show the effectiveness of the proposed strategy. 
Some conclusions and directions for future work then follow in Section~\ref{section:conclusion}.

%% file: Sections/analysis.tex
In this section, we formulate the density control problem both as a static OCP, which yields a time-independent control field that drives the swarm to a target 
equilibrium density, and as a modified dynamic OCP, which yields a time-varying control field that drives the swarm to this density within a specified time interval.
After proving an existence result for both the static and dynamic OCPs, we derive optimality conditions and prove some useful properties of the resulting system, including a global asymptotic stability result. 
Our modeling assumptions on the derivation of the macroscopic dynamics follow the same reasoning as in, e.g., \cite{sinigaglia2021density} and \cite{karthik_bil}, which we briefly review below.

We consider a swarm of robots, labeled $i = 1,...,N$, that move in a bounded domain $\Omega \in \mathbb{R}^2$. Robot $i$ occupies position $\mathbf{X}_i(t) \in \Omega$ at time $t$ and moves with controlled velocity $\mathbf{u}(\mathbf{x},t) \in \mathbb{R}^2$. This motion is perturbed by a two-dimensional Wiener process $\mathbf{W}(t)$, which models stochasticity arising from inherent sensor and actuator noise or intentionally programmed ``diffusive'' exploratory behaviors. The robot's position evolves according to the following Stochastic Differential Equation: 
\begin{equation*}
\begin{cases}
d \mathbf{X}_i(t) & = ~\mathbf{u}(\mathbf{X}_i,t)dt +\sqrt{2 \mu} d \mathbf{W}(t)+\mathbf{n}(\mathbf{X}_i(t)) d \psi(t) \\ \mathbf{X}_i(0) & = ~\mathbf{X}_{i,0},
\end{cases}
\end{equation*}
where $\mu>0$ is a diffusion coefficient, $\mathbf{n}(\mathbf{x})$ is the unit normal to the domain boundary at $\mathbf{x} \in \partial \Omega$, and $\psi(t) \in \mathbb{R}$ is a reflecting function, which ensures that the
swarm does not exit the domain. The associated probability density $q(\mathbf{x},t)$ satisfies the following linear parabolic PDE:
\begin{equation*}
    \frac{\partial q}{\partial t} + \nabla \cdot (-\mu \nabla q + \mathbf{u}q) = 0
\end{equation*}
complemented with no-flux boundary conditions. 
When considering a time-independent control field $\bar{\mathbf{u}}(\mathbf{x})$, the associated equilibrium density $\bar{q}$ satisfies
\begin{equation*}
    \nabla \cdot (-\mu \nabla \bar{q} + \bar{\mathbf{u}}\bar{q}) = 0,
\end{equation*}
which is a homogeneous linear elliptic 
advection-diffusion PDE. Note that an integral constraint of the form $\int_{\Omega} \bar{q} d\Omega = 1$ has to be added to ensure that $q$ represents a probability density, thus obtaining a well-posed problem with a nontrivial solution for each control action $\bar{\mathbf{u}}$. We can now formulate a static OCP as

\begin{equation}
\label{static_ocp}
\begin{aligned}
    &  J = \frac{\alpha}{2} \int_{\Omega} (\bar{q}-z)^2 d\Omega  + \frac{\beta}{2} \int_{\Omega} \norm{\bar{\mathbf{u}}}^2 d\Omega \quad \longrightarrow \quad \min_{\bar{q},\bar{\bb{u}}}  \\
    & \quad \quad \quad \quad  s.t.  \\
    &\begin{array}{ll}
\displaystyle \nabla \cdot (-\mu \nabla \bar{q} + \bar{\mathbf{u}}\bar{q} ) = 0 &  \textrm{in} \quad \Omega \\
\vspace{-3mm} \\
(-\mu \nabla \bar{q} + \bar{\mathbf{u}}\bar{q} )\cdot \mathbf{n} = 0 & \textrm{on} \quad \partial \Omega \\
\vspace{3mm}
\int_{\Omega} \bar{q} \, d\Omega =  1, & 
\end{array}
\end{aligned}
\end{equation}
where $\alpha,\beta >0$ are control weighting constants; $z \in L^2(\Omega)$ is the target density, which is chosen such that $\int_{\Omega} z \, d\Omega = 1$; and $\bar{q} \in H^1(\Omega)$ is the equilibrium density, which constitutes the state of our OCP. In \eqref{static_ocp}, $\bar{\mathbf{u}}\in H^1(\Omega)^2 \cap L^{\infty}(\Omega)^2$ denotes the control field, which acts bilinearly on the state dynamics as an advection field. The choice of the functional spaces will be justified in the next section. OCPs with integral state constraints are difficult to analyze and solve in general; however, in our case we can eliminate the constraint by formulating the problem in suitable zero-mean functional spaces to enforce the mass constraint effectively. 

\subsection{Analysis and functional setting}
From here on, 
we will not use the overbar to denote static variables when it is clear from the context. We briefly review some key 
properties of
problems from \cite{droniou2009noncoercive} that we can easily adapt to our case to prove asymptotic stability of the obtained optimal controls. We define the (infinite-dimensional) family of subspaces of fixed-mean functions as $\mathcal{M}_c=\{ v \in H^1(\Omega) : \int_{\Omega} v d\Omega =c\} \subset H^1(\Omega)$.
The weak formulation associated with the state problem \eqref{static_ocp} is: find $q \in \mathcal{M}_1$ such that
\begin{equation*}
    a(q,v;\mathbf{u}) = 0 \quad \forall v \in H^1(\Omega),
\end{equation*}
where the bilinear form $a$ is defined as $a(q,v;\mathbf{u}) =\int_{\Omega} \left(\mu \nabla q \cdot \nabla v - \mathbf{u}\cdot \nabla v\,q \,\right) d\Omega$. By restricting our search for $q$ to 
the space $\mathcal{M}_1$, we obtain a well-posed problem. Indeed, the state solution belongs to the kernel of the operator $L_{\mathbf{u}}:H^1(\Omega) \mapsto H^{1}(\Omega)^{*}$ defined by
\begin{equation*}
    \langle L_{\mathbf{u}} q, v \rangle = a(q,v;\mathbf{u}),
\end{equation*}
restricted to $\mathcal{M}_1 \subset H^1(\Omega)$.
In \cite{droniou2009noncoercive}, it is proven that the kernel is one-dimensional and defined up to a multiplicative constant; as a consequence, the solution is unique in $\mathcal{M}_1$ for every control velocity field $\mathbf{u}$. Furthermore, it follows from the analysis in \cite{droniou2009noncoercive} that $q>0$ a.e. on $\Omega$ and that the eigenvalues of the operator $L_{\mathbf{u}}$ are discrete and nonnegative, with the zero eigenvalue occurring with multiplicity one. Therefore, the eigenvalues of $L_{\mathbf{u}}$ when restricted to any $\mathcal{M}_c$, $c>0$, are strictly positive. We will exploit this property to prove two stability theorems for the infinite-dimensional problem and its FEM discretized counterpart.

The OCP formulation \eqref{static_ocp} does not include any weights on the spatial gradients of the control field. 
We will introduce a weight on these gradients in the cost functional such that the OCP yields a control field without steep gradients, in order to prevent control inputs whose variations are too large to be implemented on real robots. Note also that our stochastic single-integrator model will not accurately represent the microscopic dynamics of an individual robot 
if the control signal varies too quickly.
Thus, we define an auxiliary regularized problem with identical dynamics and the following cost functional:
\begin{equation*}
    J_r = J + \frac{\beta_g}{2}\int_{\Omega} ||\nabla \mathbf{u} ||^2 \, d\Omega,
\end{equation*}
where $\beta_g >0$ is the weight associated with the control gradients and $||\nabla \mathbf{u} ||$ is the Frobenius norm of $\nabla \mathbf{u}$, defined for every $\mathbf{x} \in \Omega$ as $||\nabla \mathbf{u}(\mathbf{x}) || = \sqrt{\displaystyle \sum_{i=1}^{2}\sum_{j=1}^{2} \frac{\partial u_i(\mathbf{x})}{\partial x_j}\frac{\partial u_i(\mathbf{x})}{\partial x_j}}$. 


\begin{remark}[Equivalent formulation of the state dynamics]
An equivalent way of proving the well-posedness of the state problem \eqref{static_ocp} is to make the change of variables $q = \mathring{q}+\frac{1}{|\Omega|}$, where $|\Omega|$ is the measure of $\Omega$. In this way, the integral constraint can be rewritten in terms of the new state variable $\mathring{q}$ as 
$\int_{\Omega} \mathring{q}d\Omega = 0$, a homogeneous constraint. The dynamics of $\mathring{q} \in H^1(\Omega)$ are therefore:
\begin{equation}
\label{state_dyn_eq}
\begin{array}{ll}
\displaystyle \nabla \cdot (-\mu \nabla \mathring{q} + \mathbf{u}\mathring{q} ) = - \frac{1}{|\Omega|}\nabla \cdot \mathbf{u} &  \textrm{in} \quad \Omega, \\
\vspace{-3mm} \\
\displaystyle (-\mu \nabla \mathring{q} + \mathbf{u}\mathring{q} )\cdot \mathbf{n} = -\frac{1}{|\Omega|}\,\mathbf{u}\cdot \mathbf{n} & \textrm{on} \quad \partial \Omega, \\
\vspace{1mm}
\int_{\Omega} \mathring{q} \, d\Omega =  0.
\end{array}
\end{equation}
The integral constraint can be embedded into the weak formulation by selecting $\mathcal{M}_0$ as the functional space for both test and trial functions.
The weak formulation of \eqref{state_dyn_eq} then becomes:
find $\mathring{q} \in \mathcal{M}_0$ such that \begin{equation*}
     \int_{\Omega} \mu \nabla \mathring{q} \cdot \nabla v - \mathbf{u} \cdot \nabla v \, \mathring{q} = \int_{\Omega} \frac{1}{|\Omega|}\mathbf{u} \cdot \nabla v \quad \forall v \in \mathcal{M}_0,
\end{equation*}
which we will prove to be well-posed by applying Ne\v{c}as' theorem \cite[Theorem 6.6]{salsa2015partial}. In order to do so, it is useful to note that $\norm{\nabla v}_{L^2(\Omega)}$ is a norm in $\mathcal{M}_0$ thanks to the generalized Poincar\'e inequality.
\end{remark}



In order to prove the existence of solutions to the OCP \eqref{static_ocp} with gradient regularization, it is convenient to reformulate the state dynamics with a state that  belongs to $\mathcal{M}_0$, which is a closed subspace of $H^1(\Omega)$. For any $q \in \mathcal{M}_1$, we can write the decomposition $q = w + \frac{1}{|\Omega|}$, where $w \in \mathcal{M}_0$ since $\int_{\Omega} w \,d\Omega = \int_{\Omega} q \,d\Omega-  \frac{1}{|\Omega|}\int_{\Omega}d\Omega  = 0$. In terms of $w$, the weak formulation of the problem reads: find $w \in \mathcal{M}_0$ such that
\begin{equation*}
a(w,v;\mathbf{u}) = -a(\frac{1}{|\Omega|},v;\mathbf{u}) \quad \forall v \in \mathcal{M}_0, 
\end{equation*}
where $-a(\frac{1}{|\Omega|},v;\mathbf{u}) = \int_{\Omega} \frac{1}{|\Omega|} \mathbf{u} \cdot \nabla v\,d\Omega$. The regularized cost functional weights the $H^1(\Omega)^2$-norm of the control field $\mathbf{u}$, and therefore it is natural to define the space $\mathcal{U}$ of all controls $\mathbf{u}$ as 
$\mathcal{U} = H^1(\Omega)^2 \cap L^{\infty}(\Omega)^2$.  Associated with each control $\mathbf{u}$, we can define a linear functional whose action is $F_{\mathbf{u}} v = -a(\frac{1}{|\Omega|},v;\mathbf{u}) = \int_{\Omega} \frac{1}{|\Omega|} \mathbf{u} \cdot \nabla v\,d\Omega$. The generalized Poincar\'e inequality in $H^1(\Omega)$ gives
\begin{equation*}
\norm{ w - w_{\Omega}}_{L^2(\Omega)} ~\leq~ C_p \norm{\nabla w}_{L^2(\Omega)},
\end{equation*} 
where $C_p$ is the Poincaré constant of the domain $\Omega$ and $w_{\Omega} = \int_{\Omega} w \, d\Omega = 0$. Since we therefore have 
$\norm{\nabla w}_{L^2(\Omega)} \leq \norm{w}_{H^{1}(\Omega)} \leq \sqrt{1+C_p^2} \norm{\nabla w}_{L^2(\Omega)}$, we can select the norm  $\norm{w}_{\mathcal{M}_0} = \norm{\nabla w}_{L^2(\Omega)}$. We prove that  $F_{\mathbf{u}} \in \mathcal{M}_0^{*}$, the dual of $\mathcal{M}_0$, by applying the Cauchy–Schwarz inequality and the fact that $\norm{\mathbf{u}}_{L^{2}(\Omega)^2} \leq \norm{\mathbf{u}}_{H^1(\Omega)^2}$ to obtain
\begin{equation*}
\begin{aligned}
&|F_{\mathbf{u}}v| = \left| \int_{\Omega} \frac{1}{|\Omega|} \mathbf{u} \cdot \nabla v \,d\Omega \right| \leq \frac{1}{|\Omega|}  \norm{\mathbf{u}}_{L^{2}(\Omega)^2} \norm{\nabla v}_{L^{2}(\Omega)} \\
& ~~~~~~ \leq \frac{\,\norm{\mathbf{u}}_{H^1(\Omega)^2}}{|\Omega|}   \norm{ v}_{\mathcal{M}_0},
\end{aligned}
\end{equation*} 
which implies that $\norm{F_{\mathbf{u}}}_{\mathcal{M}_0^{*}}  \leq  \frac{\norm{\mathbf{u}}_{H^1(\Omega)^2}}{|\Omega|}$.

We can cast the state equation as the following abstract variational problem: find $w \in \mathcal{M}_0$ such that
\begin{equation}
\label{var_problem}
a(w,v;\mathbf{u})  = F_{\mathbf{u}}v \quad \forall v \in \mathcal{M}_0.
\end{equation}
We will now prove the well-posedness of the variational problem \eqref{var_problem} by showing that it satisfies the hypotheses of Ne\v{c}as' theorem \cite[Theorem 6.6]{salsa2015partial}, that is, continuity and weak coercivity of the bilinear form on the left-hand side and continuity of the linear functional on the right-hand side. The most difficult property to show is the weak coercivity of the bilinear form $a$, which we prove in the following proposition. Weak coercivity is proven with respect to the pair $(\mathcal{M}_0,L_{*}^2(\Omega))$, where $L_{*}^2(\Omega)$ denotes the space of $L^2$ functions with zero mean. Note that with this choice of spaces, $\mathcal{M}_0$ is continuously and densely embedded in $L_{*}^2(\Omega)$, so that  $(\mathcal{M}_0,L_{*}^2(\Omega),\mathcal{M}_0^{*})$ is a Hilbert triplet.
\begin{proposition}[Weak coercivity of $a$]
\label{prop_1H1}
For every control $\mathbf{u} \in H^1(\Omega)^2 \cap L^{\infty}(\Omega)^2$, the bilinear form $a$ is $(\mathcal{M}_0,L_{*}^2(\Omega))$-weakly coercive, that is, there exist $\lambda > 0 $ and $\alpha > 0 $ such that 
\begin{equation*}
a(v,v;\mathbf{u}) + \lambda \int_{\Omega} v^2 \, d\Omega ~\geq~ \alpha \norm{v}_{\mathcal{M}_0}^2,
\end{equation*} 
and we can choose $\alpha = \frac{\mu}{2}$ 
and $\lambda = \frac{2}{\mu^3}  C_i^2C^4\,\norm{\mathbf{u}}^4_{H^1(\Omega)^2}$, where $C_i$ and $C$ are constants that are 
defined in the proof.
\begin{proof}
Using Holder's inequality with $(p,q,r) = (4,2,4)$, the elementary inequality $ab \leq \epsilon a^2 + \frac{b^2}{4\epsilon}$ for any $\epsilon >0$ and $a,b >0$, and the Gagliardo-Nirenberg interpolation inequality \cite[Chapter 9]{brezis}
$\norm{v}_{L^4(\Omega)}^2 \leq C_i \norm{v}_{L^2(\Omega)} \norm{\nabla v}_{L^2(\Omega)}$, where $C_i$ is the interpolation constant, we have
\begin{equation*}
\begin{aligned}
&\int_{\Omega} | \mathbf{u} \cdot \nabla v \, v | \, d\Omega ~\leq~ \norm{\mathbf{u}}_{L^4(\Omega)^2}\norm{\nabla v}_{L^2(\Omega)} \norm{v}_{L^4(\Omega)} \\
&~~~~\leq~ \epsilon \norm{\nabla v}_{L^2(\Omega)} ^2 + \frac{\norm{\mathbf{u}}^2_{L^4(\Omega)^2} \norm{v}_{L^4(\Omega)}^2}{4 \epsilon} \\ 
&~~~~\leq~ \epsilon \norm{\nabla v}_{L^2(\Omega)} ^2 + \frac{C_i \norm{\mathbf{u}}^2_{L^4(\Omega)^2} \norm{v}_{L^2(\Omega)} \, \norm{\nabla v}_{L^2(\Omega)}}{4 \epsilon} \\
&~~~~\leq~ (\epsilon + \eta)  \norm{\nabla v}_{L^2(\Omega)} ^2 + \frac{C_i^2 \norm{\mathbf{u}}^4_{L^4(\Omega)^2}}{64 \epsilon^2 \eta}\norm{v}_{L^2(\Omega)}^2 \\
&~~~~\leq~ (\epsilon + \eta)  \norm{\nabla v}_{L^2(\Omega)} ^2 + \frac{C_i^2C^4\,\norm{\mathbf{u}}^4_{H^1(\Omega)^2}}{64 \epsilon^2 \eta}\norm{v}_{L^2(\Omega)}^2
\end{aligned}
\end{equation*}
for every $\epsilon,\eta >0$, where $C>0 $ is the continuity constant of the embedding of $H^1(\Omega)^2$ into $L^4(\Omega)^2 $. Hence, it follows 
that
\begin{equation*}
\begin{aligned}
a(v,v; \mathbf{u}) + \lambda \int_{\Omega} v^2\,d\Omega &~\geq~ (\mu - \epsilon-\eta) \norm{\nabla v}^2_{L^2(\Omega)} \\ 
& \hspace{-3mm} + \left(\lambda -\frac{C_i^2C^4\,\norm{\mathbf{u}}^4_{H^1(\Omega)^2}}{64 \epsilon^2 \eta}\right) \norm{v}_{L^2(\Omega)}^2.
\end{aligned}
\end{equation*}
By setting 
$\epsilon = \frac{\mu}{4}$ and $\eta = \frac{\mu}{4}$ 
in the inequality above,
it is sufficient to choose $\lambda \geq  \frac{1}{\mu^3} C_i^2C^4\,\norm{\mathbf{u}}^4_{H^1(\Omega)^2}$ to ensure that
\begin{equation*}
\begin{aligned}
&a(v,v; \mathbf{u}) + \lambda \int_{\Omega} v^2 \, d\Omega \\
&~~~~~\geq~ \frac{\mu}{2} \norm{\nabla v}^2_{L^2(\Omega)} + \left(\lambda - \frac{C_i^2C^4\,\norm{\mathbf{u}}^4_{H^1(\Omega)^2}}{\mu^3}\right) \norm{v}_{L^2(\Omega)}^2 \\
&~~~~~\geq~ \frac{\mu}{2} \norm{\nabla v}^2_{L^2(\Omega)}  = \frac{\mu}{2} \norm{ v}^2_{\mathcal{M}_0}.
\end{aligned}
\end{equation*}
\end{proof}
\end{proposition}

We can now use Ne\v{c}as' theorem to prove the well-posedness of the state dynamics in the following theorem, which 
also provides a stability estimate that 
will be used to prove 
the existence of optimal controls.
\begin{theorem}[Existence, Uniqueness, and Stability Estimates]
\label{well_state}
For every $\mathbf{u} \in H^{1}(\Omega)^2 \cap L^{\infty}(\Omega)^2$, there exists a unique weak solution $w \in \mathcal{M}_0$ to the variational problem \eqref{var_problem} and the following stability estimate holds:
\begin{equation*}
\norm{w}_{\mathcal{M}_0} ~\leq~ \frac{2 \norm{\mathbf{u}}_{H^{1}(\Omega)^2}}{\mu \, |\Omega|}.
\end{equation*}
\begin{proof}
We verify that the hypotheses of  Ne\v{c}as' theorem are satisfied. The bilinear form $a$ is continuous since 
\begin{equation*}
\begin{aligned}
& |a(w,v)| = \left| \int_{\Omega} \Big(\mu \nabla w \cdot \nabla v - w\,\mathbf{u} \cdot \nabla v\Big) d\Omega \right| \\
&\leq~ \mu \norm{w}_{\mathcal{M}_0}\norm{v}_{\mathcal{M}_0} + \norm{\mathbf{u}}_{L^4(\Omega)^2} \norm{w}_{L^4(\Omega)} \norm{\nabla v}_{L^2(\Omega)} \\ 
& \leq~  \mu \norm{w}_{\mathcal{M}_0}\norm{v}_{\mathcal{M}_0} + C^2 \norm{\mathbf{u}}_{H^1(\Omega)^2} \norm{w}_{H^1(\Omega)} \norm{\nabla v}_{L^2(\Omega)} \\
& \leq~ \mu \norm{w}_{\mathcal{M}_0}\norm{v}_{\mathcal{M}_0} \\
& \hspace{0.5cm} + C^2 \sqrt{1+C_p^2} \norm{\mathbf{u}}_{H^1(\Omega)^2}  \norm{\nabla w}_{L^2(\Omega)} \norm{\nabla v}_{L^2(\Omega)} \\
& = ~ \Big( \mu + C^2 \sqrt{1+C_p^2}  \norm{\mathbf{u}}_{H^1(\Omega)^2} \Big) \norm{w}_{\mathcal{M}_0} \norm{v}_{\mathcal{M}_0},
\end{aligned}
\end{equation*}
where we used the generalized Poincar{\'e} inequality and 
 the continuity of the embedding of $H^1(\Omega)$ into $L^4(\Omega)$.
The bilinear form $a$ is weakly coercive according to Proposition \ref{prop_1H1}, and the linear functional $F_{\mathbf{u}}$ is continuous. Then Ne\v{c}as' theorem guarantees the existence and uniqueness of a weak solution $w \in \mathcal{M}_0$ to the variational problem \eqref{var_problem}, as well as the stability estimate 
\begin{equation*}
\norm{w}_{\mathcal{M}_0} ~\leq~ \frac{      \norm{F_{\mathbf{u}}}_{\mathcal{M}_0^{*}}     }{\alpha},
\end{equation*}
where $\alpha$ is the weak coercivity constant. Then the result follows since we can select $\alpha = \frac{\mu}{2}$ by Proposition \ref{prop_1H1} and we have shown that $  \norm{F_{\mathbf{u}}}_{\mathcal{M}_0^{*}} \leq \frac{\norm{\mathbf{u}}_{H^{1}(\Omega)^2}}{|\Omega|}   $.
\end{proof}

\end{theorem}
Note that we can also recover a stability estimate for the original state variable $q \in \mathcal{M}_1$, since 
\begin{equation}
\label{est_q}
\begin{aligned}
& \norm{q}_{H^1(\Omega)} = \norm{w + \frac{1}{|\Omega|}}_{H^1(\Omega)} ~\leq~ \norm{w}_{H^1(\Omega)}+ \norm{\frac{1}{|\Omega|}}_{H^1(\Omega)} \\ 
&\hspace{1.4cm} = ~ \norm{w}_{H^1(\Omega)} + 1 ~\leq~ \sqrt{C_p^2+1} \norm{w}_{\mathcal{M}_0} + 1  \\
&\hspace{1.4cm} \leq ~ \frac{2\,\sqrt{C_p^2+1}\norm{\mathbf{u}}_{H^{1}(\Omega)^2}}{\mu \, |\Omega|}  + 1 ~:=~ M_q.
\end{aligned}
\end{equation}

In order to prove the existence of optimal controls, we first  write a decomposition of 
the target density as $z = d + \frac{1}{|\Omega|}$, where $\int_{\Omega} d \,d\Omega = 0$ and $(q-z) = (w-d)$. Note that   
$d \in L^2_*(\Omega)$. We now consider the following optimal control problem, which is equivalent to the OCP \eqref{static_ocp}: 
\begin{equation}
\label{equiv_ocp}
\begin{aligned}
    &  J = \frac{\alpha}{2} \int_{\Omega} (w-d)^2 d\Omega  + \frac{\beta}{2} \int_{\Omega} (\norm{\mathbf{u}}^2 + \norm{\nabla \mathbf{u}}^2) d\Omega  \longrightarrow~  \min_{w,\bb{u}}  \\
    & \quad \quad \quad \quad  s.t.  \\
    &\begin{array}{ll}
a(w,v;\mathbf{u})  = F_{\mathbf{u}}v &  \forall v \in \mathcal{M}_0, \\
\end{array}
\end{aligned}
\end{equation}
where we have chosen $\beta_g = \beta$ to simplify the notation.

\begin{theorem}[Existence of Optimal Controls]
\label{ex_OCP_static}
There exists at least one optimal control pair $(w,\mathbf{u}) \in \mathcal{M}_0 \times \mathcal{U}$ for the static OCP \eqref{equiv_ocp}.
\begin{proof}
We verify that the hypotheses of Theorem 9.4 in \cite{MQS} are satisfied. 
\begin{itemize}
\item $\inf_{(w,\mathbf{u}) \in \mathcal{M}_0 \times H^1(\Omega)^2} J = \mu > -\infty$, since $J = \frac{\alpha}{2} \norm{w-d}_{L^2(\Omega)}^2 + \frac{\beta}{2} \norm{\mathbf{u}}_{H^1(\Omega)^2}^2 \geq 0$. \vspace{2mm}


\item A minimizing sequence $(w_n,\mathbf{u}_n)$ is bounded in $\mathcal{M}_0 \times \mathcal{U}$. This is because a control sequence $\{\mathbf{u}_n\}$ is bounded in $\mathcal{U}$ by definition, and we can determine that the resulting state sequence $\{w_n\}$ is bounded in $\mathcal{M}_0$ from the estimate 
in Theorem \ref{well_state},  $\norm{w_n}_{\mathcal{M}_0} \leq \frac{2 \norm{\mathbf{u}_n}_{H^{1}(\Omega)^2}}{\mu |\Omega|}$. \vspace{2mm}


\item The set of feasible state-control pairs
is weakly sequentially closed in $\mathcal{M}_0 \times \mathcal{U}$, which we demonstrate as follows. Let $\{\mathbf{u}_n\}$ be a minimizing control sequence that weakly converges to $\mathbf{u}$. Then, the resulting minimizing state sequence $\{ w_n \} $ is bounded and thus weakly convergent to $w$. Since both the state and control spaces are weakly closed, we have that $(w,\mathbf{u}) \in \mathcal{M}_0 \times \mathcal{U}$.  Define the state constraint $G(w,\mathbf{u}) \in \mathcal{M}_0^{*}$ as $\langle G(w,\mathbf{u}), v \rangle  = a(w,v;\mathbf{u})  - F_{\mathbf{u}}v $. We need to show that $G(w_n,\mathbf{u}_n) \to G(w,\mathbf{u})$ in $\mathcal{M}_0^{*}$.  Using the Cauchy–Schwarz  inequality, we obtain 
\begin{equation*}
\begin{aligned}
&|F_{\mathbf{u}_n}v - F_{\mathbf{u}}v |  ~=~ \left|\int_{\Omega} \frac{1}{|\Omega|} (\mathbf{u}_n-\mathbf{u}) \cdot \nabla v\,d\Omega \right| \\
&~~~~~\leq~ \frac{1}{|\Omega|} \norm{\mathbf{u}_n - \mathbf{u}}_{L^{2}(\Omega)^2} \norm{\nabla v}_{L^2(\Omega)} \\ 
&~~~~~ =~  \frac{1}{|\Omega|} \norm{\mathbf{u}_n - \mathbf{u}}_{L^{2}(\Omega)^2} \norm{v}_{\mathcal{M}_0} \to 0 \quad \forall v \in \mathcal{M}_0,
\end{aligned}
\end{equation*}
since $\norm{\mathbf{u}_n - \mathbf{u}}_{L^{2}(\Omega)^2} \to 0$ strongly due to the compactness of the embedding of $H^1(\Omega)^2$ into $L^2(\Omega)^2$; see, e.g., \cite[Appendix A.5.11]{MQS}. It is left to prove that
\begin{equation*}
\int_{\Omega} - \mathbf{u}_n w_n \cdot \nabla v d\Omega \to \int_{\Omega} - \mathbf{u} w \cdot \nabla v \,d\Omega \quad \forall v \in \mathcal{M}_0,
\end{equation*}
which is equivalent to showing that 
\begin{equation*}
\begin{aligned}
\int_{\Omega} (\mathbf{u}_n-\mathbf{u}) w_n \cdot \nabla v \,d\Omega & + \int_{\Omega} (w_n- w) \mathbf{u} \cdot \nabla v \,d\Omega \\
& \to 0 \quad \forall v \in \mathcal{M}_0.
\end{aligned}
\end{equation*}
For the first term, we can use the compact embedding of $H^1(\Omega)^2$ into $L^{4}(\Omega)^2$, Holder's inequality with $(p,q,r) = (4,4,2)$,  and the generalized Poincar{\'e} inequality to demonstrate that 
\begin{equation*}
\begin{aligned}
& \left|\int_{\Omega} (\mathbf{u}_n-\mathbf{u}) w_n \cdot \nabla v d\Omega  \right| \\
& ~~\leq~ \norm{\mathbf{u}_n-\mathbf{u}}_{L^{4}(\Omega)^2} \norm{w_n}_{L^{4}(\Omega)}\norm{\nabla v}_{L^{2}(\Omega)} \\
& ~~\leq~ C \,\norm{\mathbf{u}_n-\mathbf{u}}_{L^{4}(\Omega)^2} \norm{w_n}_{H^1(\Omega)}\norm{v}_{\mathcal{M}_0} \\
& ~~ \leq~  C \sqrt{1+C_p^2} \,\norm{\mathbf{u}_n-\mathbf{u}}_{L^{4}(\Omega)^2} \norm{w_n}_{\mathcal{M}_0}\norm{v}_{\mathcal{M}_0} \to 0 \quad \\
& ~~~~~~~ \forall v \in \mathcal{M}_0,
\end{aligned}
\end{equation*}
since $\norm{\mathbf{u}_n-\mathbf{u}}_{L^4(\Omega)^2} \to 0$ strongly and $ \norm{w_n}_{\mathcal{M}_0}$,  $\norm{v}_{\mathcal{M}_0}$ are bounded.
For the second term, define $\phi_v w = \int_{\Omega}  \mathbf{u} \cdot \nabla v  \, w \,d\Omega$ for every $v \in \mathcal{M}_0$ and note that $\phi_v$ is a linear and continuous functional on $\mathcal{M}_0$, since
\begin{equation*}
\begin{aligned}
|\phi_v w | & ~\leq~ \norm{\mathbf{u}}_{L^{4}(\Omega)^2} \norm{\nabla v}_{L^{2}(\Omega)} \norm{w}_{L^{4}(\Omega)} \\
& ~\leq~ C^2 \norm{\mathbf{u}}_{H^{1}(\Omega)^2} \norm{v}_{\mathcal{M}_0} \norm{ w}_{H^1(\Omega)} \\ 
& ~\leq~ C^2 \sqrt{1+C_p^2}  \norm{\mathbf{u}}_{H^{1}(\Omega)^2} \norm{v}_{\mathcal{M}_0} \norm{ w}_{\mathcal{M}_0}. 
\end{aligned}
\end{equation*}
Consequently,  $\phi_v \in \mathcal{M}_0^{*}$, and therefore we can write
\begin{equation*}
\int_{\Omega} (w_n- w) \mathbf{u} \cdot \nabla v \, d\Omega = \phi_v w_n - \phi_v w \to 0 \quad \forall v \in \mathcal{M}_0
\end{equation*}
using the definition of weak convergence in $\mathcal{M}_0$. \vspace{2mm}

\item $J$ is sequentially weakly lower semicontinuous. This can be shown by observing that since  $J$ can be written as a sum of norms, $J=\frac{\alpha}{2} \norm{w-d}_{L^2(\Omega)}^2 + \frac{\beta}{2} \norm{\mathbf{u}}_{H^1(\Omega)^2}^2$, 
it is convex and continuous on $\mathcal{M}_0 \times \mathcal{U}$ and thus 
is weakly lower semicontinuous; see, e.g., \cite[Proposition 9.1]{MQS}. \vspace{2mm}

\end{itemize}
We can now apply \cite[Theorem 9.4]{MQS}, and the result follows. 

\end{proof}
\end{theorem}

Before deriving a system of first-order necessary optimality conditions, we need to show that the control-to-state map $\Xi : \mathcal{U} \to \mathcal{M}_0$, which associates to a control $\mathbf{u}$ the resulting state $w = \Xi[\mathbf{u}]$, is Fr\'echet differentiable, that is, $\Xi \in C^1( \mathcal{U},\mathcal{M}_0)$. We prove this property in the following proposition. 

\begin{proposition}[Differentiability of the control-to-state map]
\label{z_prop}
The control-to-state map $\Xi : \mathcal{U} \to \mathcal{M}_0$ is Fr\'echet differentiable. Furthermore, its action $s = \Xi'[\mathbf{u}]\mathbf{h}$ in the control direction $\mathbf{h}$ at a state-control pair $(w,\mathbf{u})$ satisfies the following sensitivity equations:

\begin{equation}
\label{z_sens}
\begin{array}{ll}
\displaystyle-\mu \Delta s + \nabla \cdot ( s \mathbf{u} ) = - \nabla \cdot ( q \mathbf{h}) &  \textrm{in} \quad \Omega \\
\vspace{-3mm} \\
(-\mu \nabla s + \mathbf{u} s )\cdot \mathbf{n} = - q\, \mathbf{h} \cdot \mathbf{n} & \textrm{on} \quad \partial \Omega \\
\vspace{3mm}
\int_{\Omega} s \, d\Omega =  0, & 
\end{array}
\end{equation}
where $q = w + \frac{1}{|\Omega|}$ belongs to $\mathcal{M}_1$. 
\begin{proof}
The weak formulation of \eqref{z_sens} is
\begin{equation*}
a(s,v;\mathbf{u})  = F_{q \mathbf{h}} v \quad \forall v \in \mathcal{M}_0,
\end{equation*}
where $F _{q \mathbf{h}} v = \int_{\Omega} q \mathbf{h} \cdot \nabla v\,d\Omega$. We find that $F _{q \mathbf{h}} \in \mathcal{M}_0^{*}$, since 
\begin{equation*}
\begin{aligned}
|F _{q \mathbf{h}} v | & ~\leq~ \norm{\mathbf{h}}_{L^{4}(\Omega)^2} \norm{q}_{L^4(\Omega)} \norm{\nabla v}_{L^2(\Omega)} \\
& ~\leq~ C^2 \norm{\mathbf{h}}_{H^1(\Omega)^2} \norm{q}_{H^1(\Omega)} \norm{ v}_{\mathcal{M}_0} \\
& ~\leq~ C^2 M_q \norm{\mathbf{h}}_{H^1(\Omega)^2}\norm{ v}_{\mathcal{M}_0}
\end{aligned}
\end{equation*}
and $\norm{F _{q \mathbf{h}} }_{\mathcal{M}_0^{*}} \leq C^2 M_q \norm{\mathbf{h}}_{H^2(\Omega)^2}$. We can use Ne\v{c}as' theorem to establish the existence and uniqueness of solutions to the sensitivity equations \eqref{z_sens} and the stability estimate
\begin{equation*}
\norm{s}_{\mathcal{M}_0} ~\leq~ \frac{\,2\,C^2\,M_q\norm{\mathbf{h}}_{H^{1}(\Omega)^2}}{\mu } .
\end{equation*}
It remains 
to prove 
that the residual $\norm{R}_{\mathcal{M}_0} := \norm{\Xi[\mathbf{u}+\mathbf{h}] - \Xi[\mathbf{u}] - \Xi'[\mathbf{u}]\mathbf{h}}_{\mathcal{M}_0}  \to 0$ faster than  $\norm{\mathbf{h}}_{H^1(\Omega)^2}$. It can be shown that $R$ satisfies the equation
\begin{equation*}
a(R,v;\mathbf{u}+\mathbf{h}) = \int_{\Omega} s \, \mathbf{h} \cdot \nabla v \,d\Omega \quad \forall v \in \mathcal{M}_0,
\end{equation*}
where the linear functional on the right-hand side is bounded; indeed, we can obtain the following upper bound:
\begin{equation*}
\begin{aligned}
&\left|\int_{\Omega} s \, \mathbf{h} \cdot \nabla v\,d\Omega \right| ~\leq~ \norm{\mathbf{h}}_{L^{4}(\Omega)^2} \norm{s}_{L^4(\Omega)}\norm{\nabla v}_{L^2(\Omega)} \\
& \hspace{2cm} \leq~ C^2 \sqrt{1+C_p^2} \norm{\mathbf{h}}_{H^{1}(\Omega)^2}\norm{s}_{\mathcal{M}_0}\norm{v}_{\mathcal{M}_0} \\
& \hspace{2cm} \leq~  \frac{\, 2 C^4 \sqrt{1+C_p^2} \,M_q\norm{\mathbf{h}}_{H^{1}(\Omega)^2}^2}{\mu }\norm{v}_{\mathcal{M}_0} .
\end{aligned}
\end{equation*}
Then the stability estimate from Ne\v{c}as' theorem gives
\begin{equation*}
\norm{R}_{\mathcal{M}_0} ~\leq~  \frac{\, 4 C^4 \sqrt{1+C_p^2} \,M_q\norm{\mathbf{h}}_{H^{1}(\Omega)^2}^2}{\mu ^2},
\end{equation*}
and thus $\frac{ \norm{R}_{\mathcal{M}_0} }{\norm{\mathbf{h}}_{H^{1}(\Omega)^2} } \to 0 $ as $\norm{\mathbf{h}}_{H^{1}(\Omega)^2} \to 0$, which implies Fr{\'e}chet differentiability. 
\end{proof}
\end{proposition}

Note that we did not use the boundedness of the control functions $\mathbf{u}$ in the proofs of Theorems \ref{well_state} and \ref{ex_OCP_static} and Proposition \ref{z_prop}, and therefore $H^1(\Omega)^2$ could be chosen as the control space as well. However, we choose the control space 
$H^1(\Omega)^2 \cap L^{\infty}(\Omega)^2$ because  bounded controls are physically realistic for our application, and  $L^{\infty}$-norms will be used in the proof of the well-posedness of the state and control problems in the dynamic case.


\subsection{Optimality Conditions}
Following a Lagrangian approach (see, e.g., \cite[Chapter 9]{MQS}), we can recover a system of first-order necessary optimality conditions for the static problem. We will cast the optimality conditions in terms of the original state variable $q$, since this formulation is more convenient to use in the numerical treatment of the OCP. We introduce the following Lagrangian functional: 
\begin{equation*}
\mathcal{L} = J - \int_{\Omega} \nabla \cdot (-\mu \nabla q + \mathbf{u}q ) \lambda_q d\Omega + \lambda_m \Big(\int_{\Omega} q d\Omega - 1\Big),
\end{equation*}
where $\lambda_q \in H^1(\Omega)$ is the multiplier function associated with the state constraint and $\lambda_m \in \mathbb{R}$ is the scalar multiplier associated with the conservation of mass constraint, which fixes 
the unique solution to the state problem \eqref{static_ocp}.
The adjoint equation can be obtained by setting the G\^ateaux derivative of the Lagrangian functional with respect to a state variation to zero. Using integration by parts twice and substituting in the no-flux boundary conditions on the state dynamics 
(see, e.g.,\cite[Section III.B]{sinigaglia2021density} for a detailed derivation), the adjoint equation reads:  
\begin{equation*}
\begin{array}{ll}
\displaystyle  - \mu \Delta \lambda_q - \mathbf{u} \cdot \nabla \lambda_q = \alpha \, (q-z) + \lambda_m &  \textrm{in} \quad \Omega \\
\vspace{-3mm} \\
\nabla \lambda_q \cdot \mathbf{n} = 0 & \textrm{on} \quad \partial \Omega \\
\end{array}
\end{equation*}
Note that the adjoint/dual problem has a pure Neumann structure, so that we can select $\lambda_q \in \mathcal{M}_0$. It is also important to note that an explicit equation for the scalar multiplier $\lambda_m$ is not given. However,  it is shown in \cite{droniou2009noncoercive} that an additional condition on the right-hand side must hold, 
\begin{equation}
\label{adj_rhs_c}
    \int_{\Omega} \Big(\alpha \, (\bar{q}-z) + \lambda_m\Big) \hat{v}_{\mathbf{u}} \,d\Omega = 0,
\end{equation}
for $\hat{v}_{\mathbf{u}}$ in the kernel of the operator $L_{\mathbf{u}}$. Note that $q = a \hat{v}_{\mathbf{u}}$, where $a \in \mathbb{R}$ is fixed by $q$ being a probability density. Hence, the right-hand side of the adjoint equation must be orthogonal to the state. This condition is preserved by the FEM discretization and will be used to compute $\lambda_m$. 
The Euler equation is obtained by setting the G\^ateaux derivative of the Lagrangian functional with respect to a vector control variation to zero. Following a procedure similar to the adjoint derivation, we obtain:
\begin{equation}
\label{eul_eq}
    -\beta_g \Delta \mathbf{u} + \beta \, \mathbf{u} + \nabla \lambda_q q = 0.
\end{equation}
The optimal control for $\beta_g = 0$ 
has interesting properties, which we prove in the following proposition.
\begin{proposition}[Structure of the optimal control]
The optimal control solution $\mathbf{u}^{\star}$ of the OCP \eqref{static_ocp} is tangent to the boundary $\partial \Omega$ of the domain $\Omega$.

\begin{proof}
The optimal control solves the Euler equation \eqref{eul_eq}. Therefore,  $\mathbf{u}^{\star} = -\frac{1}{\beta}\nabla \lambda_q^{\star} q^{\star}$, where $\lambda_q^{\star}$ solves the adjoint equation so that we have:
\begin{equation*}
\mathbf{u}^{\star} \cdot \mathbf{n} = -\frac{1}{\beta}\nabla \lambda_q^{\star} q^{\star}\cdot \mathbf{n} = 0 \quad \textrm{on} \quad \partial \Omega
\end{equation*}
due to the adjoint boundary conditions. 
\end{proof}
\end{proposition}

The boundary condition on the optimal control enables the density dynamics to avoid obstacles. This provides an advantage over density control laws 
that only depend on the target density $z$, and thus cannot have this property.

We are now ready to formulate the associated dynamic OCP on a time interval $[0,T]$, where $T$ denotes 
both the final time and the width of the  interval, without loss of generality. It is known that for sufficiently large $T$, both the infinite-dimensional formulation of
this kind of problem and, consistently, its discrete FEM formulation
exhibit the so-called {\it turnpike behavior} \cite{trelat2015turnpike,trelat2018steady}. In other words, the dynamics of the optimal state, adjoint, and control triple 
will progress through three stages:  
first, a transient stage that 
is determined by the initial conditions; then, a ``steady-state'' stage in which the optimal triple is approximately constant and asymptotically close to an equivalent static OCP; and finally, another 
transient stage that 
is determined by the terminal conditions. 

In our setting, it is natural to require the optimal solution of the dynamic problem to converge to its static counterpart. In this way, if the initial conditions are known, we can speed up the transient to the optimal equilibrium density. Thanks to the turnpike property, it is then sufficient to set the cost functional for the dynamic problem as:
\begin{equation*}
    \begin{aligned}
    J_t & = \frac{\alpha}{2} \int_0^T \int_{\Omega} (q-\bar{q}^{\star})^2 dt\,d\Omega  + \frac{\beta}{2} \int_0^T\int_{\Omega} \norm{\mathbf{u}-\bar{\mathbf{u}}^{\star}}^2 dt\,d\Omega \\
    & + \frac{\beta_g}{2} \int_0^T\int_{\Omega} \norm{\nabla(\mathbf{u}-\bar{\mathbf{u}}^{\star})}^2 dt\,d\Omega
    \end{aligned}
\end{equation*}
with state dynamics:
\begin{equation}
\label{state_dyn_t}
\begin{array}{ll}
\displaystyle \frac{\partial q}{\partial t} + \nabla \cdot (-\mu \nabla q + \mathbf{u}q ) = 0 &  \textrm{in} \quad \Omega \times (0,T) \\
\vspace{-3mm} \\
(-\mu \nabla q + \mathbf{u}q )\cdot \mathbf{n} = 0 & \textrm{on} \quad \partial \Omega \times (0,T) \\
\vspace{-3mm} \\
q(\mathbf{x},0) = q_0(\mathbf{x}) & \textrm{on} \quad \Omega \times \{0\}
\end{array}
\end{equation}
and the additional box constraint $|\mathbf{u}(\mathbf{x},t)| \leq \norm{\bar{\mathbf{u}}^{\star}}_{L^{\infty}(\Omega)^2}$ for every $\mathbf{x} \in \Omega$ and  a.e. on $t \in (0,T)$, which ensures that  the magnitude of the dynamic control action does not exceed the magnitude of  its static counterpart, in an $L^{\infty}$ sense. Note that the dynamic OCP weights the distance between the dynamic state variable $q(t)$ and its optimal static equilibrium. Because of this, the final transient due to the turnpike property is eliminated, ensuring convergence to the static optimal state-control pair.



The analysis of the dynamic OCP follows similar arguments to the ones given, 
e.g., in \cite[Proposition 6]{roy2018fokker}, where control functions are chosen in 
$L^2(0,T;H_0^1(\Omega)^2)$ with box constraints. Note that this is equivalent to selecting controls in $L^2(0,T;H_0^1(\Omega)^2 \cap L^{\infty}(\Omega)^2)$.
In our case, we do not require the controls to be zero at the domain boundaries, and so we select the 
control space $\mathcal{U}_t = L^2(0,T;\mathcal{U})$. Using the same decomposition of the density $q$ as in 
the static OCP, we consider the dynamics of $w(t) = q(t) - \frac{1}{|\Omega|}$ with initial conditions $w(\mathbf{x},0) = q(\mathbf{x},0) - \frac{1}{|\Omega|}$. Due to the no-flux boundary conditions, the set $\mathcal{M}_1$ is forward invariant for the dynamics \eqref{state_dyn_t}; that is, the system is mass-conservative. It is easy to show this property: defining $m(t) = \int_{\Omega} q(\mathbf{x},t) d\Omega$, we have that
\begin{equation}
\label{conserve}
\begin{aligned}
    \dot{m} & ~=~ \int_{\Omega} \frac{\partial q}{\partial t} d\Omega  ~=~ - \int_{\Omega}\nabla \cdot (-\mu \nabla q + \mathbf{u}q )d\Omega \\
&~=~ \int_{\partial \Omega} (-\mu \nabla q + \mathbf{u}q )\cdot \mathbf{n}d\Omega ~=~ 0.
\end{aligned}
\end{equation}
As a consequence, $w(t) \in \mathcal{M}_0$ for any $t>0$ and any choice of the dynamic control $\mathbf{u} \in \mathcal{U}_t$. Note also that the integral constraint on the density $q$ need not 
be taken into account explicitly, due to the mass-preserving property of the state equation and the fact that $q_0 \in \mathcal{M}_1$. Another fundamental property of the density dynamics \eqref{state_dyn_t} is that the state solution is nonnegative, that is, $q(t) \geq 0 $ for a.e. $t \in (0,T)$; see, e.g., \cite[Chapter 7]{evans}, for a proof that uses  maximum principles for parabolic equations.

We can now cast the dynamic OCP in terms of $w$ in the same way as the static OCP  and select the 
state space $\mathcal{Y} = H^1(0,T;\mathcal{M}_0,\mathcal{M}_0^{*})$. The weak formulation of the dynamic state problem reads: find $w \in \mathcal{Y}$ such that 
\begin{equation*}
\left\langle \frac{\partial w(t)}{\partial t},v \right\rangle + a(w(t),v;\mathbf{u}(t)) = F_{\mathbf{u}(t)} v \quad \forall v \in \mathcal{M}_0  
\end{equation*}
for a.e. $t \in (0,T)$. The following proposition establishes the well-posedness of the state dynamics and a stability estimate that 
will be used in the proof of the existence of optimal controls for the dynamic OCP. 
\begin{proposition}
For every control $\mathbf{u} \in \mathcal{U}_t$, there exists a unique solution $w \in \mathcal{Y}$ that 
satisfies the stability estimate
\begin{equation*}
\begin{aligned}
\norm{w}^2_{\mathcal{Y}} ~\leq~&~  \frac{4 M^2 + 1}{\mu}e^{2 \lambda T} \norm{w_0}^2_{L^2(\Omega)} \\ 
& ~+ \Big( \frac{8 M^2 + 1}{\mu^2}e^{2 \lambda T} +2 \Big)\frac{\norm{\mathbf{u}}^2_{L^2(0,T;H^1(\Omega)^2)}}{|\Omega|^2},
\end{aligned}
\end{equation*}
where the constants $M$ and $\lambda$ are defined in the proof.
\begin{proof}
Consider the Hilbert triplet $(\mathcal{M}_0,L_{*}^2(\Omega),\mathcal{M}_0^{*})$ and note that the hypotheses of \cite[Theorem 7.1]{MQS} are satisfied since $F_{\mathbf{u}} \in L^2(0,T;\mathcal{M}_0^*)$; indeed,
\begin{equation*}
\begin{aligned}
&\norm{F_{\mathbf{u}}}_{L^2(0,T;\mathcal{M}_0^*)}^2 ~=~ \int_{0}^{T} \norm{F_{\mathbf{u}(t)}}_{\mathcal{M}_0^{*}}^2 \, dt \\ 
& \quad \quad \quad \quad \leq~  \int_{0}^{T}  \frac{\norm{\mathbf{u}(t)}^2_{H^1(\Omega)^2}}{|\Omega|^2} ~=~ \frac{\norm{\mathbf{u}}^2_{L^2(0,T;H^1(\Omega)^2)}}{|\Omega|^2}
\end{aligned}
\end{equation*}
and $w_0 = q_0 - \frac{1}{|\Omega|} \in L^2(\Omega)$. The bilinear form $a(q,v,\mathbf{u})$ is continuous for a.e. $t \in (0,T)$ since
\begin{equation*}
\begin{aligned}
& |a(q,v;\mathbf{u}(t))|  \leq \left( \mu + \sqrt{1+C_p^2}  \norm{\mathbf{u}(t)}_{L^4(\Omega)^2} \right) \norm{w}_{\mathcal{M}_0} \norm{v}_{\mathcal{M}_0} \\
&~~~\leq  \Big( \mu + C_{\infty} \sqrt{1+C_p^2}  \norm{\mathbf{u}(t)}_{L^\infty(\Omega)^2}\Big) \norm{w}_{\mathcal{M}_0} \norm{v}_{\mathcal{M}_0} \\
&~~~\leq  \Big( \mu + C_{\infty} \sqrt{1+C_p^2}  \norm{\bar{\mathbf{u}}^{\star}}_{L^{\infty}(\Omega)^2}\Big) \norm{w}_{\mathcal{M}_0} \norm{v}_{\mathcal{M}_0} 
\end{aligned}
\end{equation*}
due to the definition of the control space $\mathcal{U}_t$, where $C_{\infty}$ is the continuity constant of the embedding of $L^{\infty}(\Omega)^2$ into $L^4(\Omega)^2$. We define $M = \left( \mu + C_{\infty} \sqrt{1+C_p^2}  \norm{\bar{\mathbf{u}}^{\star}}_{L^{\infty}(\Omega)^2}\right)$ as the continuity constant of the bilinear form.
The bilinear form is also $\mathcal{M}_0 - L^2(\Omega)$ weakly coercive for a.e. $t\in (0,T)$, and we can choose the weak coercivity constant independently of 
time as $\lambda = \frac{2}{\mu^3} C_i^2 C_{\infty} ^4 \norm{\bar{\mathbf{u}}^{\star}}_{L^{\infty}(\Omega)^2}^4$, which can be obtained from Proposition \ref{prop_1H1} using the embedding of $L^{\infty}(\Omega)^2$ into $L^{4}(\Omega)^2$ and the $L^{\infty}$ box constraint on the dynamic control in terms of its static counterpart. It is also clear that the bilinear form is $t$-measurable, so that the stability estimate follows by applying \cite[Theorem 7.1]{MQS}  in our setting.

\end{proof}
\end{proposition}

Using 
similar arguments as for the static problem,  it is possible to show that the control-to-state map is also differentiable 
in the dynamic problem  and that there exists at least one optimal control. Since the proof is similar to the proof of Theorem \ref{ex_OCP_static}, we just sketch the 
part that is particular to the dynamic problem  in the proof of the following proposition.
\begin{proposition} \label{prop:ExistOptCtrlPair}
There exists at least one optimal control pair $(w,\mathbf{u}) \in \mathcal{Y} \times \mathcal{U}_t$ for the dynamic OCP \eqref{state_dyn_t} written in terms of the zero-mean function $w$.
\begin{proof}(Sketch)
In order to apply \cite[Theorem 9.4]{MQS}, we need to show that that the set of feasible points is weakly sequentially closed in $\mathcal{Y} \times \mathcal{U}_t$. This requires proving the nontrivial result that 
\begin{equation*}
\int_{0}^{T} \int_{\Omega} \mathbf{u}_n q_n \cdot \nabla v \, d\Omega \, dt ~\to~ \int_{0}^{T} \int_{\Omega} \mathbf{u} q \cdot \nabla v \, d\Omega \, dt  
\end{equation*}
for every $ v \in L^2(0,T,\mathcal{M}_0)$, which is equivalent to proving 
that
\begin{equation*}
\begin{aligned}
& \int_{0}^{T} \int_{\Omega} (\mathbf{u}_n-\mathbf{u}) w_n \cdot \nabla v\,d\Omega\,dt \\ 
& \hspace{2cm} + ~ \int_{0}^{T} \int_{\Omega} (w-w_n) \mathbf{u} \cdot \nabla v \,d\Omega\,dt  ~ ~\to ~~ 0 
\end{aligned}
\end{equation*}
 for every $ v \in L^2(0,T,\mathcal{M}_0)$.
For the first term, 
define $f_n(t) = \int_{\Omega} (\mathbf{u}_n(t)-\mathbf{u}(t)) w_n(t) \cdot \nabla v(t)\,d\Omega$ for a.e.  $t \in (0,T)$ and for every $v \in \mathcal{M}_0$. By the definition of the control and state spaces, we have that $\mathbf{u}_n(t),\mathbf{u}(t) \in H^1(\Omega)^2$ and $w_n(t) \in \mathcal{M}_0$. Using the same reasoning as for the static problem, we find that 
\begin{equation*}
|f_n(t)| ~\leq~ M(t) \norm{\mathbf{u}_n(t)-\mathbf{u}(t)}_{L^4(\Omega)^2}  ~\to~ 0,
\end{equation*}
where $M(t) = C \sqrt{1+C_p^2} \norm{w_n(t)}_{\mathcal{M}_0}  \norm{v(t)}_{\mathcal{M}_0}$ is bounded for a.e. $t$. Then, by the Lebesgue dominated convergence theorem, we obtain 
\begin{equation*}
\int_{0}^{T} \int_{\Omega} (\mathbf{u}_n-\mathbf{u}) w_n \cdot \nabla v  \, d\Omega dt ~\to~ 0 \quad  \forall v \in L^2(0,T,H^1(\Omega)).
\end{equation*}
For the second term, we have that
\begin{equation*}
\begin{aligned}
&\left|\int_{0}^{T} \int_{\Omega} (w-w_n) \mathbf{u} \cdot \nabla v \,d\Omega\,dt \right| \leq \\
& \norm{\bar{\mathbf{u}}^{\star}}_{L^{\infty}(\Omega)^2} \norm{w-w_n}_{L^2(0,T;L^2(\Omega)}\norm{\nabla v}_{L^2(0,T;L^2(\Omega)} ~\to~ 0,
\end{aligned}
\end{equation*}
since $\norm{\bar{\mathbf{u}}^{\star}}_{L^{\infty}(\Omega)^2}$ and $ \norm{\nabla v}_{L^2(0,T;L^2(\Omega)}$ are bounded and $\norm{w-w_n}_{L^2(0,T;L^2(\Omega))} \to 0 $. Indeed, the space $\mathcal{Y}$ is compactly embedded in $L^2(0,T;L_{*}^2(\Omega))$ according to Aubin's Lemma; see, e.g., \cite[Theorem A.19]{MQS}. 
\end{proof} 
\end{proposition}


The adjoint equation for this modified dynamic problem is:
\begin{equation*}
\begin{array}{ll}
\displaystyle -\frac{\partial \lambda_q}{\partial t} - \mu \Delta \lambda_q - \mathbf{u} \cdot \nabla \lambda_q = \alpha \, (q-\bar{q}^{\star})  &  \textrm{in} \quad \Omega \times (0,T) \\
\vspace{-3mm} \\
\nabla \lambda_q \cdot \mathbf{n} = 0 & \textrm{on} \quad \partial \Omega \times (0,T) \\
\lambda_q(\mathbf{x},T) = 0 & \textrm{on} \quad \Omega \times \{T\},
\end{array}
\end{equation*}
while the Euler equation or reduced gradient can be written as
\begin{equation}
    \nabla J = -\beta_g \Delta (\mathbf{u}-\bar{\mathbf{u}}^{\star}) + \beta \, (\mathbf{u}-\bar{\mathbf{u}}^{\star}) + \nabla \lambda_q q.
\end{equation}

Moreover, for a sufficiently large time interval $(0,T)$, the turnpike property implies that after a transient due to the initial conditions, $\mathbf{u}^{\star}(t) \to \bar{\mathbf{u}}^{\star}$ and $q^{\star}(t) \to \bar{q}^{\star}$. This can be interpreted as a robustness property: even if the initial conditions are not exactly known, the convergence to the optimal steady-state control action still ensures convergence to the target density, since the equilibrium induced by $\bar{\mathbf{u}}^{\star}$ is globally asymptotically stable. The global asymptotic stability of this equilibrium is established by the following theorem.

\begin{theorem}
\label{stab_cont}
The optimal solution  $\bar{q}^{\star}(\bar{\mathbf{u}}^{\star})$ of problem \eqref{static_ocp} is globally asymptotically stable on the mass-preserving subspace $\mathcal{M}_1.$

\begin{proof}
Existence and uniqueness of $\bar{q}^{\star}(\bar{\mathbf{u}}^{\star})$ follows from Theorem \ref{well_state}. To show global asymptotic stability, consider the Lyapunov function $\mathcal{V}=\frac{1}{2}\int_{\Omega} (q-\bar{q}^{\star})^2 d\Omega$. Define the error $e:=q-\bar{q}^{\star}$ and note that $e \in \mathcal{M}_0$. $\mathcal{V}$ is positive on $\mathcal{M}_0$ and vanishes only for $q=\bar{q}^{\star}$. The time derivative of $\mathcal{V}$ along the solution of the state equation is:
\begin{equation*}
\begin{aligned}
    \dot{\mathcal{V}} &=   \int_{\Omega} e \, \frac{\partial q}{\partial t} d\Omega =- \int_{\Omega} e \, \nabla \cdot (-\mu \nabla q+\bar{\mathbf{u}}^{\star}q) d\Omega \\
    &= -\int_{\Omega} e \, \nabla \cdot (-\mu \nabla e+\bar{\mathbf{u}}^{\star}e) d\Omega \\
    &=- \int_{\Omega}  (\mu \norm{\nabla e}^2 - \bar{\mathbf{u}}^{\star}\cdot \nabla e\,e ) \, d\Omega = -\langle L_{\bar{\mathbf{u}}^{\star}} e, e \rangle,
\end{aligned}
\end{equation*}
where we have substituted the steady-state condition $\nabla \cdot (-\mu \bar{q}^{\star} + \bar{\mathbf{u}}^{\star}\bar{q}) = 0 $ and the no-flux boundary conditions. Now, the final expression of $\dot{\mathcal{V}}$ is a weak formulation associated with the state operator, and it is formulated on $\mathcal{M}_0$. It is proven in \cite{droniou2009noncoercive} that this operator is strictly positive on this subspace, and thus $\dot{\mathcal{V}}<0$.
\end{proof}
\end{theorem}

In the next section, we will prove that the FEM discretization inherits the properties of the associated infinite-dimensional formulation, and thus a stability proof for the semi-discrete system can use arguments similar to Theorem \ref{stab_cont} but in a finite-dimensional, algebraic setting.

%% file: Sections/numerical_analysis.tex
In this section, we analyse the FEM discretization of both the static and dynamic OCPs and propose numerical algorithms for their solutions that 
are 
based on the algebraic properties of the finite-dimensional problem.
The FEM discretization inherits the structure of the infinite-dimensional problem and, in particular, the state dynamics reduce to a kernel-finding problem with a unique solution determined by the discretized integral mass constraint.

In the static problem, the FEM discretization of the state dynamics is
\begin{equation}
\label{state_fem}
\begin{array}{l}
\Big(A - \mathbb{B}_x^{\top}\mathbf{u}_x -\mathbb{B}_y^{\top}\mathbf{u}_y \Big)\mathbf{q}=\mathbf{0}, \\
\mathbf{F}^{\top}\mathbf{q} - 1 = \mathbf{0}; \\
\end{array}
\end{equation}
the adjoint dynamics are given by 
\begin{equation}
\label{adj_fem}
\begin{array}{l}
\Big(A - \mathbb{B}_x\mathbf{u}_x -\mathbb{B}_y\mathbf{u}_y \Big)\bs{\lambda}_q =\alpha M (\mathbf{q}-\mathbf{z}) + \lambda_m \mathbf{F}; \\
\end{array}
\end{equation}
and the Euler equation can be expressed as
\begin{equation*}
\begin{aligned}
    & \beta M_u \mathbf{u}_x + \beta_g A_u \mathbf{u}_x + \bs{\lambda}_q^{\top}\mathbb{B}_x \mathbf{q} = \mathbf{0}, \\
    &\beta M_u \mathbf{u}_y + \beta_g A_u \mathbf{u}_y +  \bs{\lambda}_q^{\top}\mathbb{B}_y \mathbf{q} = \mathbf{0}.
\end{aligned}
\end{equation*}
In these equations, $A$ and $M$ are the usual stiffness and mass matrices, while 
$\mathbb{B}_x$ is a rank-3 tensor defined as $\mathbb{B}_{x,ijk} = \int_{\Omega} \frac{\partial \phi_i}{\partial x} \phi_j \phi_k \, d\Omega$ and $\mathbb{B}_y$ is defined in a similar way.
For an alternative and equivalent way of defining the FEM matrices, see \cite{sinigaglia2021density}. The vectors $\bb{q}$, $\bb{z}$, $\bs{\lambda}_q \in \mathbb{R}^{N_q}$ denote the coefficients of the FEM basis functions for the state, target, and adjoint variables, respectively. The entries of the vector $\mathbf{F} \in \mathbb{R}^{N_q}$ are defined as $\mathbf{F}_i = \int_{\Omega} \phi_i \, d\Omega$,  where $\phi_i$ is the corresponding FEM basis function. Given this definition, the mass of an FEM variable $v_h = \sum_{i=1}^{N} \phi_i v_i$ is simply:
\begin{equation*}
    \int_{\Omega} v_h d \Omega = \sum_{i=1}^{N} \int_{\Omega} \phi_i \, d\Omega \,v_i = \mathbf{F}^{\top} \mathbf{v}.
\end{equation*}
With a slight abuse of notation, we group the FEM discretizations of the $x$ and $y$ components of the control action as $\mathbf{u} = [\mathbf{u}_x \,\, \mathbf{u}_y]^{\top}$ and define $\mathbb{B}=\texttt{Stack3}(\mathbb{B}_x,\mathbb{B}_y)$, where the $\texttt{Stack3}$ operation stacks $\mathbb{B}_x$ and $\mathbb{B}_y$ along the third direction. In this way, we can compactly write
\begin{equation*}
    \mathbb{B}_x \mathbf{u}_x + \mathbb{B}_y \mathbf{u}_y =  \mathbb{B} \mathbf{u} ;
\end{equation*}
note that $\mathbb{B} \mathbf{u} $ is a matrix. The transpose operation on $\mathbb{B}$ is defined as $\mathbb{B}_{ijk}^{\top} = \mathbb{B}_{jik}$.
The transpose of the state matrix,  $A-\mathbb{B}\mathbf{u}$, can be easily seen to be the discrete adjoint matrix, $A$ being symmetric. This, in turn, implies full commutativity of the Discretize-then-Optimize (DtO) and Optimize-then-Discretize (OtD) solution methods for this problem; see \cite{sinigaglia2021density} and \cite[Ch. 6]{MQS}, for more details.

The semi-discrete set of modified optimality conditions arising from the discretization (in space) of the dynamic problem is: 

\begin{subequations}
    \begin{eqnarray}
    \label{q_t}
    M \dot{\q} + (A-\mathbb{B}^{\top}\mathbf{u})\q = \mathbf{0}, &  t \in (0,T) \nonumber \\
    \hspace{0.7cm} \q(0) = \q_0; 
    \end{eqnarray} 
    \begin{eqnarray}
    \label{l_t}
    -M \dot{\bs{\lambda}}_q + (A-\mathbb{B}\mathbf{u})\bs{\lambda}_q =     \alpha M_q(\q-\bar{\q}^{\star})  , &  t \in (0,T) \nonumber \\
    \hspace{0.7cm} \bs{\lambda}_q(T) = \mathbf{0};
    \end{eqnarray} 
    \begin{eqnarray}
    \Big(\beta M_u +\beta_g A_u\Big) (\mathbf{u}_x-\bar{\mathbf{u}}_x^{\star}) + \bs{\lambda}_q^{\top} \mathbb{B}_x^{\top} \q = \mathbf{0} , &  t \in (0,T) \nonumber \\
    \Big(\beta M_u +\beta_g A_u\Big) (\mathbf{u}_y-\bar{\mathbf{u}}_y^{\star}) + \bs{\lambda}_q^{\top} \mathbb{B}_y^{\top} \q = \mathbf{0} , &  t \in (0,T) \nonumber \\ \notag
    \end{eqnarray}
\end{subequations}
where we have used the definition of $\mathbb{B}$ to compactly write both the state and adjoint dynamics. We remark that in the static problem, we solve for a constant vector $\bar{\mathbf{u}} \in \mathbb{R}^{2N_u}$, while in the dynamic problem, our unknowns $\mathbf{q}=\mathbf{q}(t),\mathbf{u}=\mathbf{u}(t)$ are time-dependent.

\subsection{Properties of FEM discretization}

We define the mass-preserving linear subspace arising from the FEM discretization as 
$\tilde{\mathcal{M}}_c = \{ \mathbf{v} \in \mathbb{R}^{N_q} : \mathbf{F}^{\top} \mathbf{v} = c \}$ for some total mass $c > 0$. In the following, we will need in particular $\tilde{\mathcal{M}}_0$ and $\tilde{\mathcal{M}}_1$; note the close parallel 
with their infinite-dimensional counterparts $\mathcal{M}_0$ and $\mathcal{M}_1$. We can now prove a number of useful properties that the FEM approximation inherits from the infinite-dimensional problem, thus making it a consistent (and elegant) discretization of the OCP. 
\begin{proposition}
\label{adj_one}
For each $\mathbf{u} \in \mathbb{R}^{2N_u}$, $\mathbf{1} \in \texttt{Ker}(A - \mathbb{B} \uc )$ 
and the dimension of $\texttt{Ker}(A - \mathbb{B} \uc )$ is 1; that is, $\texttt{Span}(\texttt{Ker}(A - \mathbb{B} \uc )) = \{\mathbf{1}\}$.

\begin{proof}
The matrix $A - \mathbb{B} \uc$ 
is the FEM discretization of the adjoint PDE operator, 
which is defined up to a constant 
since the adjoint system is a pure Neumann problem. In FEM terms, this constant corresponds to the vector $\mathbf{1} \in \mathbb{R}^{N_q}$; see, e.g., \cite{Bochev2005}. 
\end{proof}
\end{proposition}

From Proposition \ref{adj_one}, a simple yet useful result follows for the discretized state problem, also ensuring its well-posedness. This result is stated in the following proposition. 

\begin{proposition}
\label{prop_v}
The kernel of the state matrix $(A-\mathbb{B}^{\top}\mathbf{u})$ is one-dimensional, that is, $\texttt{Dim}(\texttt{Ker}(A-\mathbb{B}^{\top}\uc))=1$  
$\forall \mathbf{u} \in \mathbb{R}^{2N_u}$.

\begin{proof}
$\texttt{Rank}(A - \mathbb{B}^{\top}\mathbf{u}) = \texttt{Rank}(A - \mathbb{B}\mathbf{u}) = N_q-1$, and thus $\texttt{Dim}(\texttt{Ker}(A-\mathbb{B}^{\top}\uc))=1$.
\end{proof}
\end{proposition}

We denote the vector spanning the kernel of $A-\mathbb{B}^{\top}\mathbf{u}$ as $\mathbf{v}(\mathbf{u})$ and note that the kernel of $A-\mathbb{B}\mathbf{u}$ is spanned by the vector of ones, $\mathbf{1}$, for every control action $\mathbf{u} \in \mathbb{R}^{2N_u}$. The previous results allow us to prove the following proposition regarding the mass-preserving property of the semi-discrete system.

\begin{proposition}
\label{semi-conserve}
The FEM discretization of the state equation in the dynamic OCP, 
\begin{equation*}
\begin{array}{ll}
    M \dot{\q} + (A-\mathbb{B}^{\top}\mathbf{u})\q = \mathbf{0}, & t \in (0,T) \\
    \mathbf{q}(0)=\mathbf{q}_0,
\end{array}
\end{equation*}
is mass-conservative with respect to the FEM mass function $m_d(t) = \int_{\Omega} \sum_{i=1}^{N_q} \phi_i q_i d\Omega = \mathbf{F}^{\top} \q$. That is, $\dot{m}_d = 0$ for every control action $\mathbf{u}$.

\begin{proof}
We can use the relation $\mathbf{F} = M \mathbf{1}$ (see, e.g., \cite{Bochev2005}), Proposition \ref{adj_one}, and the semi-discrete state dynamics to obtain:
\begin{equation*}
    \dot{m}_d = \mathbf{F}^{\top} \dot{\q} = \mathbf{1}^{\top} M \dot{\q} = -\mathbf{1}^{\top}(A-\mathbb{B}^{\top}\mathbf{u})\q = 0.
\end{equation*}
\end{proof}
\end{proposition}

Proposition \ref{semi-conserve} constitutes the finite-dimensional analogue of Equation \eqref{conserve}.  In other words, by choosing $\mathbf{q}_0 \in \tilde{\mathcal{M}}_1$, we have that $\mathbf{q}(t) \in \tilde{\mathcal{M}}_1$ for a.e. $t \in (0,T)$, which means that the linear subspace $\tilde{\mathcal{M}}_1$ is forward invariant for the semi-discrete dynamics. 
We are now ready to prove a useful stability theorem which ensures that for every control action $\bar{\mathbf{u}}$, the resulting equilibrium is unique and globally asymptotically stable on the relative mass-preserving subspace, which we assume to be $\tilde{\mathcal{M}}_1$. Without loss of generality, we prove this result for the equilibrium density induced by the optimal control $\bar{\mathbf{u}}^{\star}$.

\begin{theorem}
\label{stab_discrete}
For every optimal control $\bar{\mathbf{u}}^{\star}$, the resulting optimal equilibrium density $\bar{\q}^{\star}(\bar{\mathbf{u}}^{\star})$ is unique and globally asymptotically stable.

\begin{proof}
Since $\texttt{Dim}(\texttt{Ker}(A-\mathbb{B}^{\top}\bar{\uc}^{\star}))=1$, the equilibrium density has the form $\bar{\mathbf{q}}^{\star}(\bar{\mathbf{u}}^{\star}) = k \mathbf{v}(\bar{\mathbf{u}}^{\star})$, where $\mathbf{v}(\bar{\mathbf{u}}^{\star})$ is the vector spanning the one-dimensional kernel of the state dynamics and $k\in \mathbb{R}$ has to be determined. Imposing the condition $\bar{\mathbf{q}}^{\star} \in \tilde{\mathcal{M}}_1$, that is $\mathbf{F}^{\top}\bar{\mathbf{q}}^{\star}(\bar{\mathbf{u}}^{\star}) = 1$, we obtain a unique solution.

The bilinear form $a(q,v;\mathbf{u})$ that arises from the weak formulation of the state equation is associated with 
an operator with nonnegative eigenvalues when $q,v$ belong to the zero-mean space. We proved that the FEM solution $\q(t)$ remains on the  mass-preserving subspace $\tilde{\mathcal{M}}_1$. Defining $B \in \mathbb{R}^{N_q \times N_{q}-1}$ as a basis for $\tilde{\mathcal{M}}_0 \subset \mathbb{R}^{N_q}$,  there exist vectors $\dot{\mathbf{w}}(t),\mathbf{\bar{w}}^{\star} \in \mathbb{R}^{N_q-1}$ such that $\dot{\q}(t) = B \dot{\mathbf{w}}(t)$ and $\bar{\mathbf{q}}^{\star} = B \bar{\mathbf{w}}^{\star} + \frac{1}{|\Omega|}\mathbf{1}$. 
Furthermore, since the FEM approximation selects $q_h,v_h \in V_h \subset H^1(\Omega)$, the strict positivity of the operator on the mass-preserving subspace implies that:
\begin{equation*}
    a(v_h,v_h;\bar{\mathbf{u}}^{\star}) = \mathbf{v}^{\top} B^{\top}(A-\mathbb{B}^{\top} \bar{\mathbf{u}}^{\star})B \mathbf{v} > 0
\end{equation*}
for every $\mathbf{v} \in \mathbb{R}^{N_q-1}$.
Now, consider the candidate Lyapunov function $l(\q) = \frac{1}{2} (\q-\bar{\q}^{\star})^{\top} M (\q-\bar{\q}^{\star})$. It is clear that $l>0 \quad \forall \q \, \neq \bar{\q}^{\star}$, since the mass matrix $M$ is positive definite. Using the previous results, we show that $\dot{l}< 0$ along solutions of the semi-discrete FEM dynamics $M \dot{\q} + (A-\mathbb{B}^{\top} \bar{\mathbf{u}}^{\star})\q = \mathbf{0}$. Indeed, since the bilinear form $a(q,v;\mathbf{u})$ is strictly positive on $\mathcal{M}_0$,\,$\tilde{\mathcal{M}}_0 \subset \mathcal{M}_0$ and $\mathbf{q}(t)-\bar{\q}^{\star} = B \mathbf{w}(t)$ for some $\mathbf{w}(t) \in \mathbb{R}^{N_q-1}$, we have that
\begin{equation*}
\begin{aligned}
    \dot{l} 
    &= (\q-\bar{\q}^{\star})^{\top} M \dot{\q} \\
    &= -(\q-\bar{\q}^{\star})^{\top} (A-\mathbb{B}^{\top} \bar{\mathbf{u}}^{\star})\q \\
    &= -(\q-\bar{\q}^{\star})^{\top} (A-\mathbb{B}^{\top} \bar{\mathbf{u}}^{\star})(\q-\bar{\q}^{\star}) \\
    &= -\mathbf{w}^{\top}B^{\top} (A-\mathbb{B}^{\top} \bar{\mathbf{u}}^{\star}) B \mathbf{w} < 0.
\end{aligned}
\end{equation*}
\end{proof}
\end{theorem}
Note that it is not unexpected that the finite-dimensional stability proof mirrors the steps of its infinite-dimensional counterpart (Theorem \ref{stab_cont}), since the FEM discretization inherits the properties of the associated infinite-dimensional operators.

When solving the dynamic OCP, a consistent discretization in time is needed. The so-called $\theta$-method (see, e.g., \cite[Ch. 8]{MQS}) for semi-discrete evolution problems reads, in our case,
\begin{equation}
\label{theta_method}
\begin{array}{l}
 \begin{aligned}
    &\frac{M}{\Delta t}(\mathbf{q}_{i+1}-\mathbf{q}_{i}) + \theta (A-\mathbb{B}^{\top}\mathbf{u}_{i+1})\mathbf{q}_{i+1} \\
    & ~~~ + (1-\theta)(A-\mathbb{B}^{\top}\mathbf{u}_{i})\mathbf{q}_{i} = \mathbf{0},  \quad i = 0,\ldots,N_t - 1 
\end{aligned} \\
\\
 \mathbf{q}_0 \quad \textrm{given}
\end{array}
\end{equation}
where $\Delta t$ is the time step, $N_t$ is the number of time discretization points, and $\theta \in [0,1]$ is selected by the user. The choice of $\theta=0$ and $\theta=1$ gives the classical forward and backward Euler methods, respectively. It can also be proven that the only method of order $2$ in time is the Crank-Nicolson method, which is obtained by selecting $\theta = \frac{1}{2}$; see, e.g., \cite[Ch. 5]{quart}. Note that $\q_0$ is set equal to the initial condition.
A straightforward application of Proposition \ref{adj_one} allows us to prove that the full discretization using a $\theta$-method remains mass-preserving.
\begin{proposition}
\label{prop_theta}
The $\theta$-method time discretization is mass-preserving. That is, if $\q_{i} \in \tilde{\mathcal{M}}_1$, then  $\q_{i+1} \in \tilde{\mathcal{M}}_1$ as well, or equivalently,
\begin{equation*}
    \mathbf{F}^{\top}\big(\q_{i+1}-\q_{i}\big) = 0.
\end{equation*}
\begin{proof}
Since $\mathbf{1}$ is in the kernel of $(A-\mathbb{B}^{\top}\mathbf{u})$ for every $\mathbf{u}\in \mathbb{R}^{2 N_u}$, we can use the same arguments as in Proposition \ref{semi-conserve} to obtain:
\begin{equation*}
\begin{aligned}
    & \mathbf{F}^{\top}\big(\q_{i+1}-\q_{i}\big) = \mathbf{1}^{\top} M \big(\q_{i+1}-\q_{i}\big) \\
    & ~~~ = -\Delta t \,\theta\, \mathbf{1}^{\top} (A-\mathbb{B}^{\top}\mathbf{u}_{i+1})\mathbf{q}_{i+1} \\
    & ~~~~~~ -\Delta t \,(1-\theta)\, \mathbf{1}^{\top} (A-\mathbb{B}^{\top}\mathbf{u}_{i})\mathbf{q}_{i} = 0. 
\end{aligned}
\end{equation*}
\end{proof}
\end{proposition}
In other words, Proposition \ref{prop_theta} states that the increment is orthogonal to the FEM mass vector $\mathbf{F}$, and thus the discrete evolution of $\q_i$ remains in the mass-preserving subspace for all $i=0,\ldots,N_t$ and for every control action $\mathbf{u}_i$, which may be time-varying.

The last step of our numerical analysis is to ensure that the numerical solution $\mathbf{q}_i$ remains nonnegative for every $i=1,\ldots,N_t$, given that $\q_0\geq \mathbf{0}$. In order to rigorously prove this property {\it a priori} for a fully discrete scheme, a discrete form of the maximum principle needs to be proved. Positivity is ensured if the stiffness matrix resulting from the FEM discretization is an M-matrix (see, e.g., \cite{refId0,xu1999monotone}) and the mass lumping technique is used together with, for example, a backward Euler method \cite{thomee2014positivity,xu1999monotone}. Modified finite-element schemes are available to ensure these properties under mild assumptions on the geometry of the finite element mesh, see, e.g., \cite{xu1999monotone,zhang2013maximum}. However, for sufficiently fine meshes of strict Delaunay type, the standard FEM together with the mass lumping technique is positivity preserving (see \cite{xu1999monotone,zhang2013maximum,thomee2008existence}). For the sake of simplicity, here we use a standard Galerkin method with mass lumping, making sure that the FEM mesh satisfies the geometric assumptions that ensure positivity preservation. We remark that in each of our numerical test cases, the density remains nonnegative  throughout all optimization iterations. 

\subsection{Solution algorithm}
Exploiting the properties of the algebraic systems governing the state and adjoint variables, we can derive a numerical algorithm to compute the reduced gradient. The main difficulty is to find the Lagrange multiplier $\lambda_m$ associated with 
the mass constraint. By projecting the discrete adjoint equation on the kernel of the state equation, we can recover an equation for $\lambda_m$. Then, some care is needed in the numerical treatment of the adjoint system. Since it results from the discretization of a pure Neumann problem, the adjoint problem comprises a singular system with a one-dimensional kernel spanned by $\mathbf{1}$ \cite{Bochev2005}. This, in turn, means that its solution is defined up to an arbitrary constant. Following \cite{Bochev2005}, a more robust way to solve the adjoint system is to look for solutions which have zero mean. In the infinite-dimensional formulation,  this amounts to requiring that $\int_{\Omega} \lambda_q d\Omega = 0$, which in the FEM discretization readily translates to $\mathbf{F}^{\top}\bs{\lambda}_q=0$. Note the duality with the state system, which should satisfy $\mathbf{F}^{\top}\mathbf{q}=1$. For the sake of brevity, we just provide the FEM discretization of the optimization problem whose solution provides the adjoint system. An associated variational formulation in continuous space can also be derived (see, e.g., \cite{Bochev2005}). For fixed controls, the adjoint solution solves the linearly constrained, quadratic optimization problem:
\begin{equation}
\label{adj_opt}
\begin{aligned}
    &  \,\, \frac{1}{2}\,\bs{\lambda}_q^{\top} \left(A - \mathbb{B} \uc  \right)\bs{\lambda}_q - \Big(\alpha\, M_q (\q-\mathbf{z})+ \lambda_m \mathbf{F}\Big)^{\top} \bs{\lambda}_q  \longrightarrow \min_{\bs{\lambda}_q}\\
    & \vspace{0.5cm} \\
    & s.t. \, \begin{array}{l}
\mathbf{F}^{\top} \bs{\lambda}_q = 0, \\
\end{array}
\end{aligned}
\end{equation}
where $\lambda_m$ can be computed using the kernel properties of the state system and the adjoint system. The KKT system arising from Problem \eqref{adj_opt} is a well-posed sparse linear system. The numerical computation of the reduced gradient is sketched in Algorithm \ref{alg_red_grad}.

\begin{algorithm}[t]
  \begin{algorithmic}[1]
\State $\mathbf{v}(\mathbf{u}) \gets \texttt{Span}\Big(\texttt{Ker}\left( A - \mathbb{B}^{\top} \uc \right)\Big)$  \Comment{Solve state equation \eqref{state_fem}}
\State $\q \gets \frac{\mathbf{v}(\mathbf{u})}{\mathbf{F}^{\top}\mathbf{v}(\mathbf{u})}$  
\State $\lambda_m \gets -\frac{\alpha\, \mathbf{v}(\uc)^{\top}M_q (\q-\mathbf{z})}{\mathbf{v}(\uc)^{\top}\mathbf{F}}$ 




\State $\bs{\lambda}_q \gets  \textrm{Solve Problem \eqref{adj_opt}}$


\State $\nabla_{\uc_x}\tilde{J} = \beta M_u \mathbf{u}_x + \bs{\lambda}_q^{\top} \mathbb{B}_x^{\top} \q$  \Comment{Compute reduced gradient}
\State $\nabla_{\uc_y}\tilde{J} = \beta M_u \mathbf{u}_y + \bs{\lambda}_q^{\top} \mathbb{B}_y^{\top} \q$ 

  \end{algorithmic}
  \caption{Reduced gradient computation with integral density constraint }
  \label{alg_red_grad}
\end{algorithm}

The explicit computation of $\lambda_m$ in Algorithm \ref{alg_red_grad} utilizes the kernel properties of the state matrix. Left-multiplying the adjoint equation \eqref{adj_fem} by $\mathbf{v}(\mathbf{u})$, and applying the result in Proposition \ref{prop_v}, we obtain
\begin{equation*}
 \alpha \mathbf{v}(\mathbf{u})^{\top}M (\mathbf{q}-\mathbf{z}) + \lambda_m \mathbf{v}(\mathbf{u})^{\top}\mathbf{F} =  \mathbf{v}(\mathbf{u})^{\top}\Big(A - \mathbb{B}\mathbf{u} \Big)\bs{\lambda}_q = 0,
\end{equation*}
which can be used to compute $\lambda_m$. Note that this condition ensures the well-posedness of the adjoint equation and corresponds to Equation \eqref{adj_rhs_c}.

\begin{figure}
    \centering
    \includegraphics[width= \linewidth]{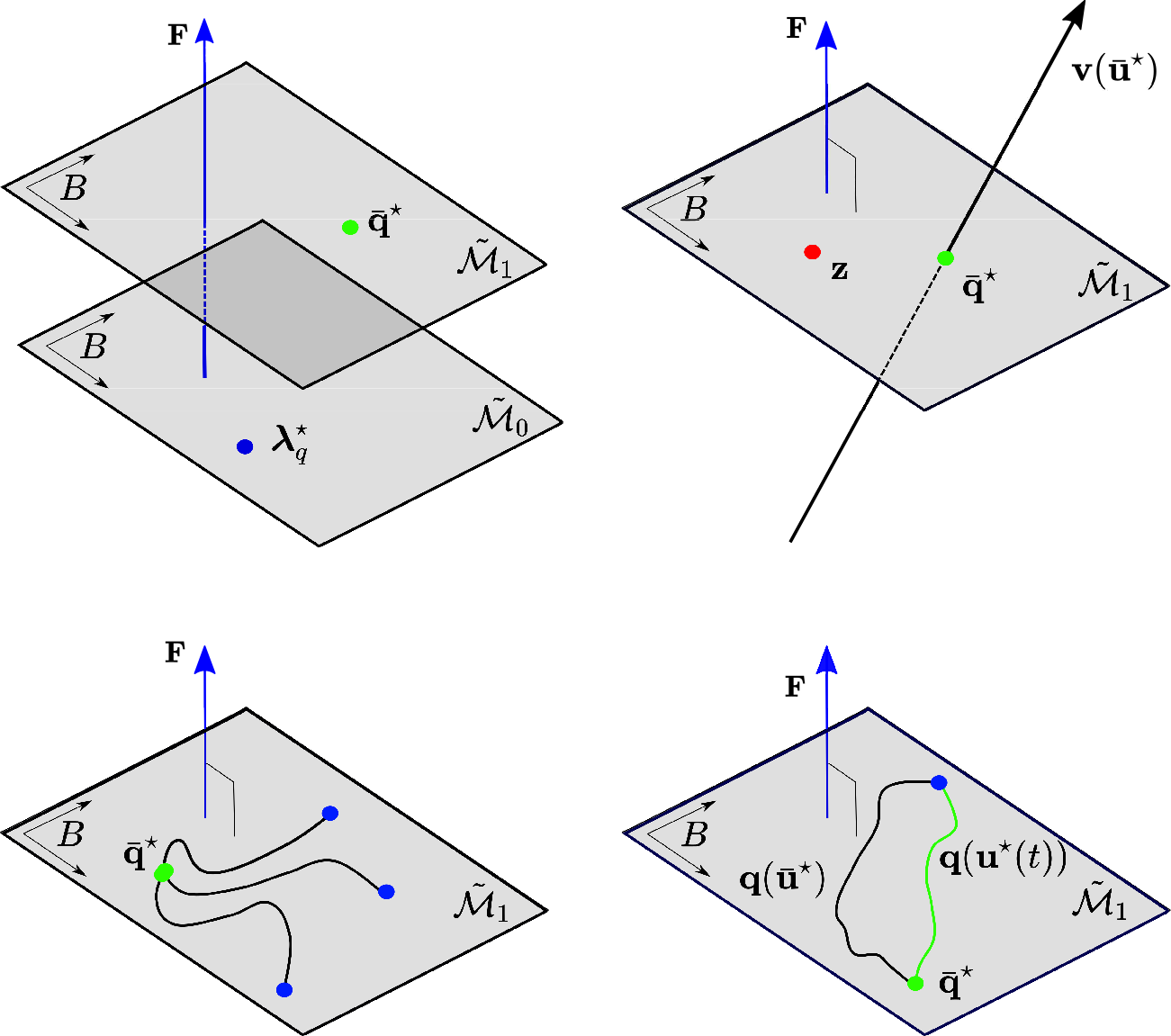}
    \caption{Illustration of the main properties of the FEM discretization. {\it Top left:} Duality between state and adjoint when the adjoint is solved with a zero-mean constraint. {\it Top right:} Interpretation of the state problem as an eigenvector optimization problem for the kernel of the state operator, which is spanned by $\mathbf{v}(\mathbf{u})$. {\it Bottom left:} From any initial condition in $\tilde{\mathcal{M}}_1$, the state asymptotically converges to $\bar{\mathbf{q}}^{\star}$ while remaining in $\tilde{\mathcal{M}}_1$. {\it Bottom right:} Dynamic optimization procedure exploiting knowledge of the initial condition and of the static optimal solution.}
    \label{fig:sketch}
\end{figure}

In order to compute the solution to the static optimization problem, we use the quasi-Newton method outlined in Algorithm \ref{alg_newton}. In this method, 
the reduced Hessian is approximated by a matrix of the form $\beta M_u + \beta_g A_u$. Note that this matrix is not the FEM equivalent of 
the reduced Hessian, which is both 
harder to derive and 
more difficult to compute. 
For a derivation of second-order necessary conditions for a similar problem, see \cite{aronna2021first}; an alternative numerical treatment based on the conjugate gradient method, which does not require the computation of the reduced Hessian, can 
be found in \cite{glowinski2021bilinear}. The dynamic problem is solved with a similar iterative quasi-Newton method, whose steps are sketched in Algorithm \ref{alg_newton_t}. 
We illustrate the properties of the discretized OCP in Figure \ref{fig:sketch}. 

\begin{algorithm}[t]
  \begin{algorithmic}[1]
  \State $H \gets \beta M_u + \beta_g A_u$
    \For{$i=0:\text{maxIter}$} 
        \State $\nabla J(\bb{u}^{(i)}) \gets 
        \textrm{Algorithm1}(\bb{u}^{(i)})$ \Comment{Compute reduced gradient}
        \State $\bb{d}^{(i)} \gets \text{Solve } H \bb{d}^{(i)} = -\nabla J(\bb{u}^{(i)})  $ 
        \State $\tau \gets \text{ArmijoBacktracking}(J,\bb{d}^{(i)},\bb{u}^{(i)}) $ \Comment{Line search}
        \State $\bb{u}^{(i+1)} \gets \bb{u}^{(i)} + \tau \bb{d}^{(i)}  $ \Comment{Update control}
        \If{ $\norm{\nabla J(\bb{u}^{(i)})} < \text{tol}$}
        \State $  \text{return}$
        \EndIf
    \EndFor
  \end{algorithmic}
  \caption{Modified Newton method for static OCP}
  \label{alg_newton}
\end{algorithm}
 
\begin{algorithm}[t]
  \begin{algorithmic}[1]
  \State $H \gets \beta M_u + \beta_g A_u$
  \State $\bar{\mathbf{q}}^{\star},\bar{\mathbf{u}}^{\star},\bar{\bs{\lambda}}^{\star} \gets \textrm{Algorithm2}$ \Comment{Solve static OCP}
  \State $\mathbf{u}^{0} \gets \textrm{Repmat}(\bar{\mathbf{u}}^{\star},N_t)$
 \Comment{$N_t$ copies of solution to static OCP} 
    \For{$i=0:\text{maxIter}$} 
        \State $\mathbf{q}^{(i)} \gets \textrm{SolveStateDyn}(\mathbf{u}^{(i)})$ \Comment{Solve Eq. \eqref{q_t}}  
        \State $\bs{\lambda}_q^{(i)} \gets \textrm{SolveAdjointDyn}(\mathbf{u}^{(i)},\mathbf{q}^{(i)},\bar{\q}^{\star})$
        \Comment{Solve Eq. \eqref{l_t}}
        \State $\nabla J(\bb{u}^{(i)}) \gets 
        H (\mathbf{u}^{(i)}-\bar{\mathbf{u}}^{\star}) + \bs{\lambda}_q^{(i)^\top} \mathbb{B}^{\top} \q^{(i)}$ 

        \State $\bb{d}^{(i)} \gets \text{Solve } H \bb{d}^{(i)} = -\nabla J(\bb{u}^{(i)})  $ 
        \State $\tau \gets \text{ArmijoBacktracking}(J,\bb{d}^{(i)},\bb{u}^{(i)}) $ \Comment{Line search}
        \State $\bb{u}^{(i+1)} \gets \bb{u}^{(i)} + \tau \bb{d}^{(i)}  $ \Comment{Update control}
        \If{ $\norm{\nabla J(\bb{u}^{(i)})} < \text{tol}$}
        \State $  \text{return}$
        \EndIf
    \EndFor
  \end{algorithmic}
  \caption{Modified Newton method for dynamic OCP  }
  \label{alg_newton_t}
\end{algorithm}

%% file: Sections/numerical_simulations.tex
In this section, we show the effectiveness of our control algorithm through numerical simulations of three test cases. In all cases, the computational domain is discretized into a triangular mesh with $N_q$ degrees of freedom, and the time interval $[0,T]$ is discretized into time steps $\Delta t = 0.03$ [s],
where 
$T = 3$ [s]. 
Test case 2 is also run for $T = 100$ [s] with the control field $\bar{\mathbf{u}}^{\star}$ from the static OCP only.
The resulting fully discrete optimization problem has $N_q$ state variables and $2N_q$ control variables in the static case, while in the dynamic case the number of variables is multiplied by $N_t$. We set $N_q = 2704$ in test case 1, $N_q = 2603$ in test case 2, and $N_q = 3831$ in test case 3.
Computations are carried out in MATLAB using a modified version of the \texttt{redbKit} library \cite{quarteronireduced} to assemble the FEM matrices and tensors and the \texttt{TensorToolbox} \cite{tensor_toolbox} to perform efficient tensor computations. 

In all test cases, the diffusion coefficient is normalized to  $\mu=1$, and the control weightings $\alpha$, $\beta$, and $\beta_g$ are selected using a trial-and-error procedure to obtain satisfactory tracking performance. We note that the diffusion coefficient $\mu$ influences the $L^2$-norm of the optimal tracking error, $\int_{\Omega} (\bar{q}^{\star}-z)^2 d\Omega$: for higher $\mu$, a stronger control field is needed to constrain the optimal equilibrium density $\bar{q}^{\star}$ to a given distance from the target density $z$, requiring the control weighting $\beta$ to be reduced.

\begin{figure}[t] 
\centering
\includegraphics[width=0.5\textwidth]{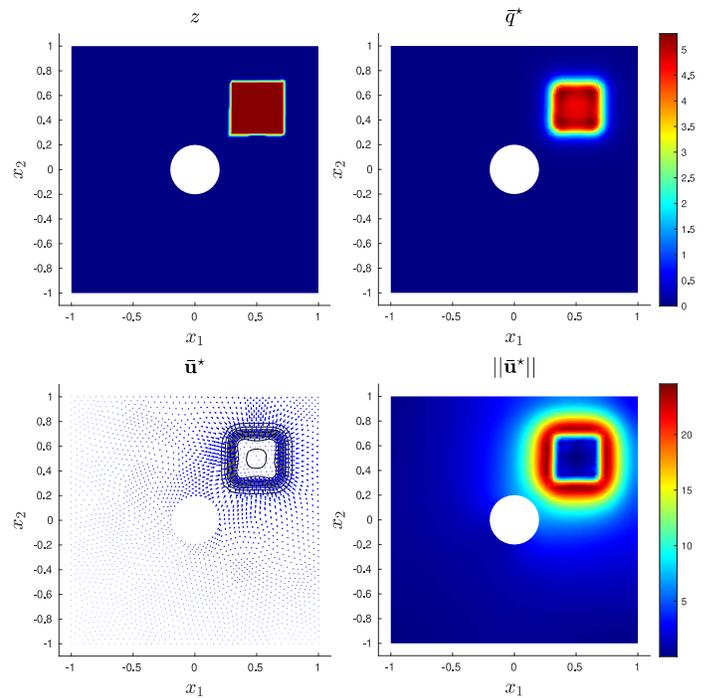}
\caption{Test case 1. Equilibrium density $\bar{q}^{\star}$ and control field $\bar{\mathbf{u}}^{\star}$ 
computed from the static optimization problem, in which $z$ is a non-smooth target function. The control weights 
are 
$\alpha=1$, $\beta=10^{-3}$, and $\beta_g=10^{-5}$.}
\label{tc1_static}
\end{figure}

\begin{figure}[h!]
\centering

\includegraphics[width=0.4\textwidth]{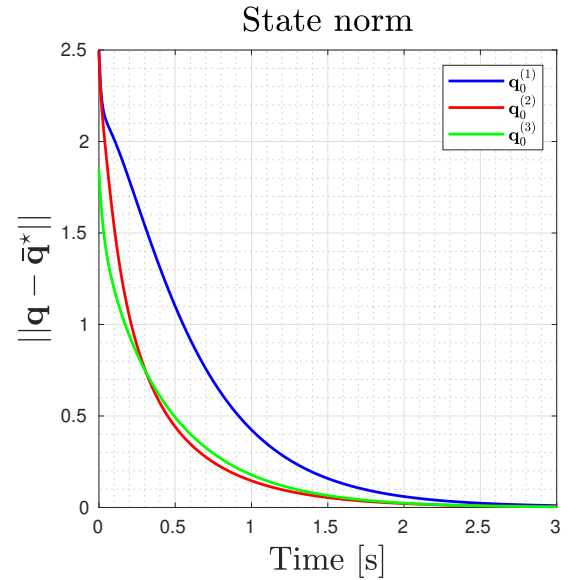}

\caption{Test case 1. Convergence in the $L^2$ norm to the optimal equilibrium $\bar{\mathbf{q}}^{\star}$ from three different initial conditions, $\mathbf{q}^{(1)}_0$, $\mathbf{q}^{(2)}_0$, and $\mathbf{q}^{(3)}_0$, where $\mathbf{q}^{(1)}_0$ and $\mathbf{q}^{(2)}_0$ are Gaussian densities centered at $(-0.5,-0.5)$ and $(-0.5,0.5)$, respectively, and $\mathbf{q}^{(3)}_0$ is the uniform density over the domain $\Omega$.}
\label{tc1_static_conv}
\end{figure}

\begin{figure}[h!]
\centering
\includegraphics[width=0.5\textwidth]{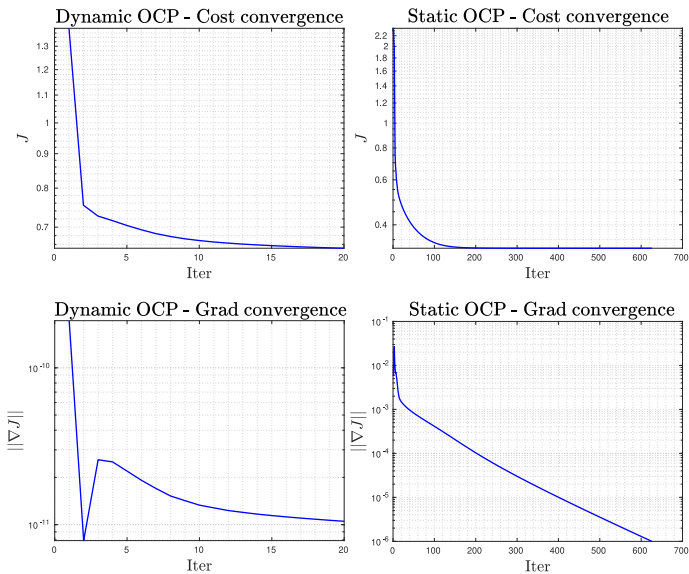}

\caption{Test case 1. Convergence of the cost functional $J$ ({\it top}) and $||\nabla J||$ ({\it bottom}) during iterations of the quasi-Newton method used to solve
the dynamic ({\it left}) and static ({\it right}) optimization problems. 
}
\label{t_sim_conv}
\end{figure}

\begin{figure}[h!]
\centering
\includegraphics[width=0.5\textwidth]{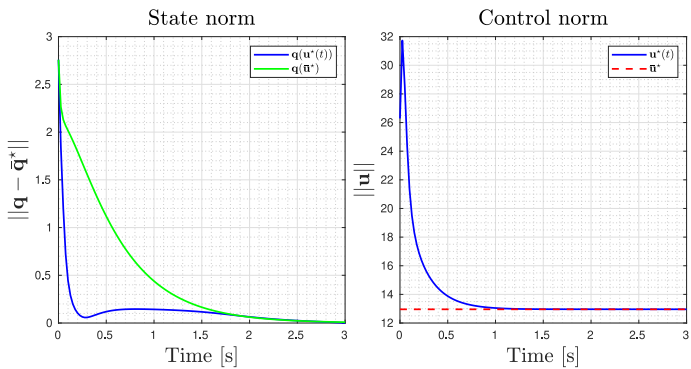}

\caption{Test case 1. Convergence in the $L^2$ norm of the solution of the dynamic optimization problem to its static counterpart. {\it Left}: Norm convergence of the density $\q$ to the optimal equilibrium $\bar{\mathbf{q}}^{\star}$ under the constant control field $\bar{\mathbf{u}}^{\star}$ (green) and the time-varying control field $\mathbf{u}^{\star}(t)$ (blue). 
{\it Right}: Norm convergence of $\mathbf{u}^{\star}(t)$ to $\bar{\mathbf{u}}^{\star}$.
}
\label{t_sim_speedup}
\end{figure}

\begin{figure}[h!]
\centering
\includegraphics[width=0.5\textwidth]{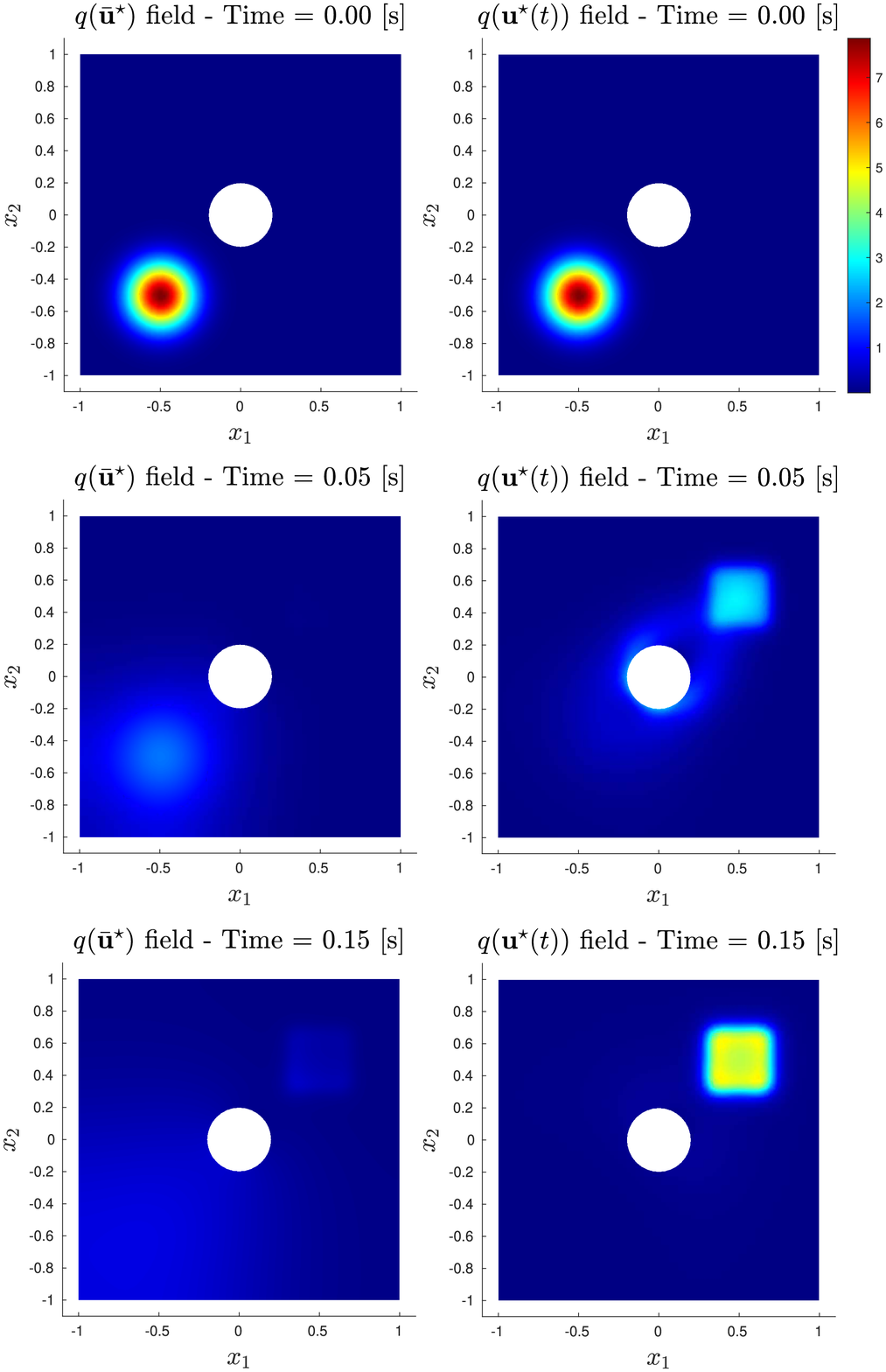}
\caption{Test case 1. Density evolution under ({\it left}) the optimal constant control field $\bar{\mathbf{u}}^{\star}$, and ({\it right}) the optimal time-varying control field $\mathbf{u}^{\star}(t)$, which is optimized for the  initial condition but converges to $\bar{\mathbf{u}}^{\star}$ over time. 
}
\label{t_sim_frames_q}
\end{figure}

\begin{figure}[h!]
\centering
\includegraphics[width=0.5\textwidth]{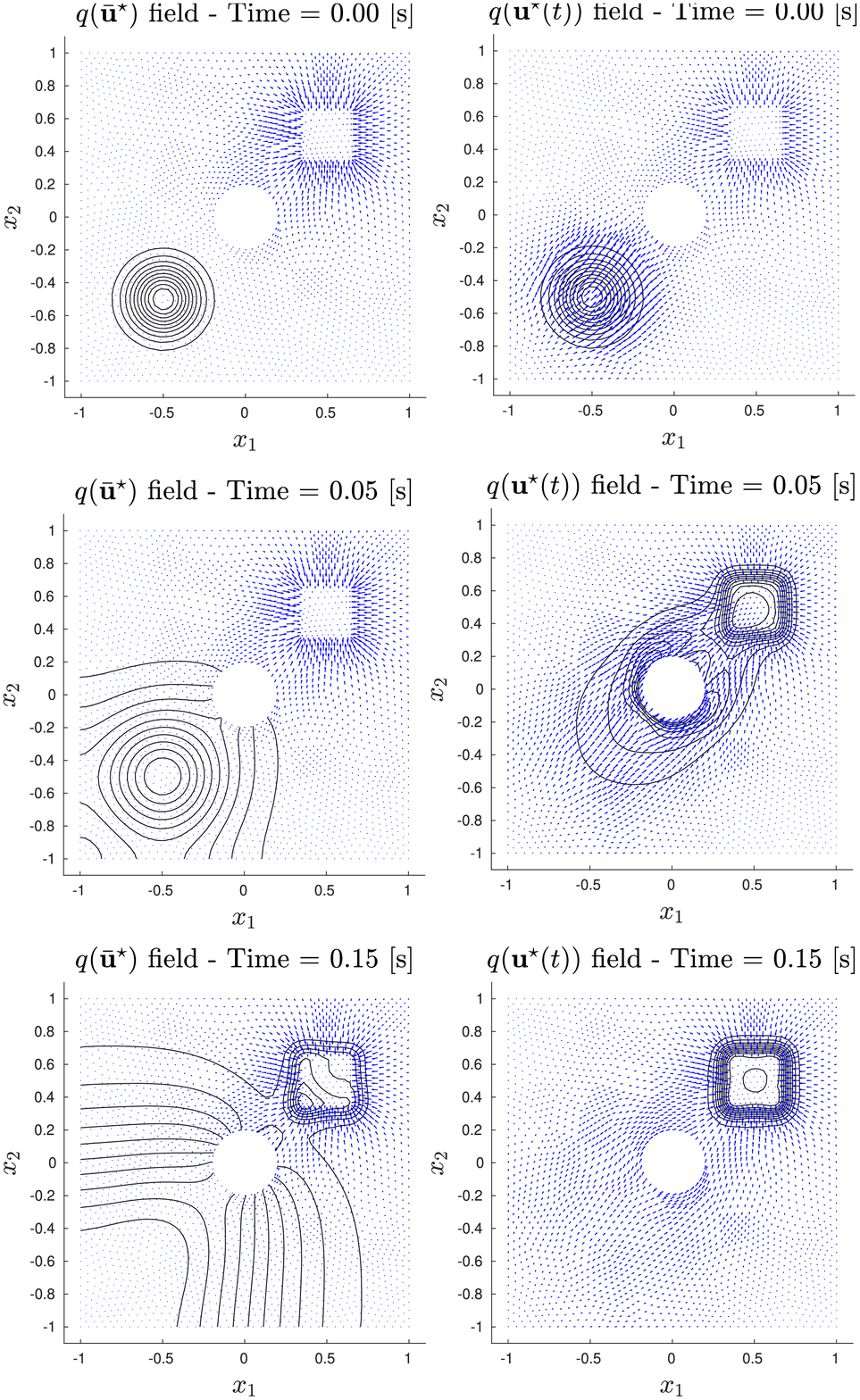}
\caption{Test case 1. Evolution of density contours (black) under the optimal control 
fields (blue), $\bar{\mathbf{u}}^{\star}$ ({\it left}) and  $\mathbf{u}^{\star}(t)$ ({\it right}).
}
\label{t_sim_frames_u}
\end{figure}

In test case 1, we define the domain as 
$\Omega=[-1,1]^2 \setminus B(\mathbf{0},0.2)$, where $B(\mathbf{x},r)$ denotes the two-dimensional ball centered at $\mathbf{x}$ with radius $r$. 
The ball $B(\mathbf{0},0.2)$ represents a circular obstacle which must be avoided by the density dynamics. 
Figure \ref{tc1_static} plots the target density $z$ and the numerical solution of the static optimization problem, comprised of the equilibrium density $\bar{q}^{\star}$ and control field $\bar{\mathbf{u}}^{\star}$, along with the norm of $\bar{\mathbf{u}}^{\star}$. The plots show that the equilibrium density is smooth, non-negative, and close to the target density, and that the control field varies smoothly over the domain and does not have steep gradients. We simulated the system under the constant control field $\bar{\mathbf{u}}^{\star}$ for three different initial densities. Figure  \ref{tc1_static_conv} shows that for each initial condition, $\bar{\mathbf{u}}^{\star}$ stabilizes the discretized state $\q$ to the corresponding optimal equilibrium density $\bar{\q}^{\star}$. 
We then solve the modified dynamic problem to optimize the convergence rate to 
equilibrium from the initial condition $\mathbf{q}_{0}^{(1)}$, which is defined as a Gaussian density centered at $(-0.5,-0.5)$. Convergence of the solutions 
of both the static and dynamic problems to the optimal solution is evidenced by 
the plots of the cost functional $J$ and $||\nabla J||$ in Figure \ref{t_sim_conv}. 
Figure \ref{t_sim_speedup} shows that the time-varying control field $\mathbf{u}^{\star}(t)$ produces faster convergence to the optimal equilibrium density than $\bar{\mathbf{u}}^{\star}$ and that 
$\mathbf{u}^{\star}(t) \to \bar{\mathbf{u}}^{\star}$ for sufficiently large $T$.
Finally, Figures \ref{t_sim_frames_q} and \ref{t_sim_frames_u} plot several 
snapshots of the density evolution under both control fields, $\bar{\mathbf{u}}^{\star}$ and $\mathbf{u}^{\star}(t)$. 
Figure \ref{t_sim_frames_u} shows that 
$\mathbf{u}^{\star}(t)$ at $t = 0$ [s] exhibits its highest magnitude near the peak of the initial density; this concentration of control effort enables it to drive the swarm around the obstacle to the target equilibrium density faster than 
$\bar{\mathbf{u}}^{\star}$. 

\begin{figure}[h!]
\centering
\includegraphics[width=0.5\textwidth]{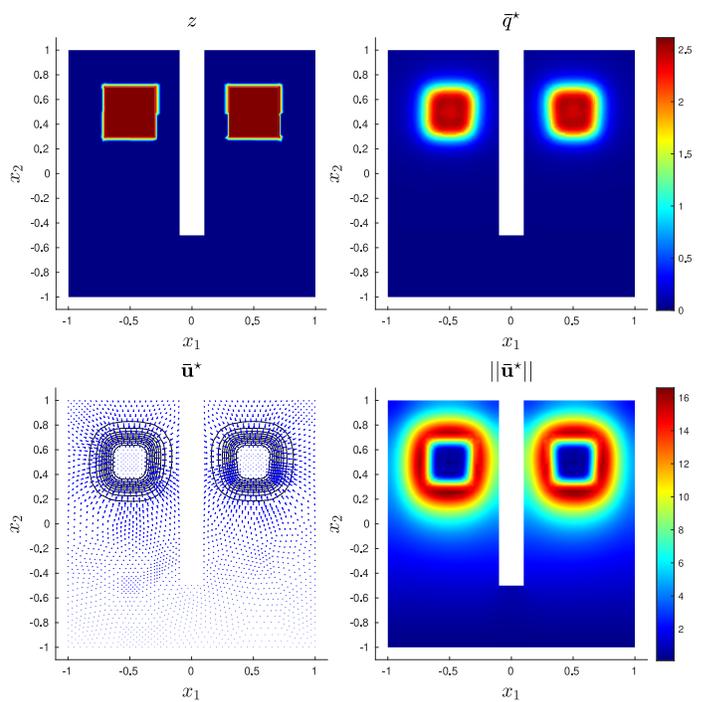}
\caption{Test case 2. Equilibrium density $\bar{q}^{\star}$ and control field $\bar{\mathbf{u}}^{\star}$ 
computed from the static optimization problem, in which $z$ is a non-smooth target function
composed of the disjoint union of two characteristic functions. The peaks of $z$ 
are separated by an obstacle, which prevents density transport directly between the peaks. 
The control weights 
are 
$\alpha=1$, $\beta=10^{-3}$, and $\beta_g=10^{-5}$.}
\label{fig:tc2_static}
\end{figure}

\begin{figure}[h!]
\centering
\includegraphics[width=0.5\textwidth]{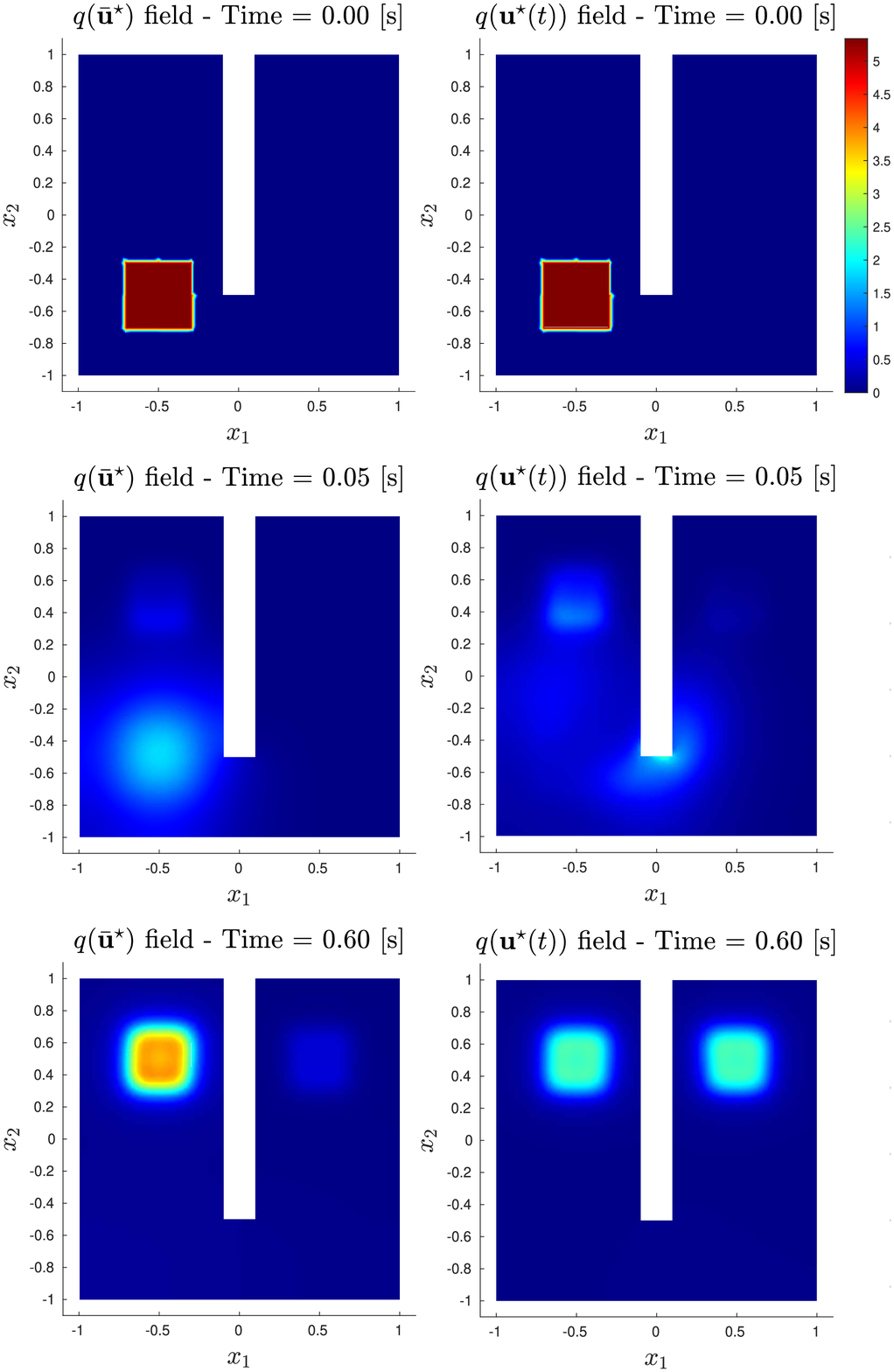}
\caption{Test case 2. Density evolution under ({\it left}) the optimal constant control field $\bar{\mathbf{u}}^{\star}$, and ({\it right}) the optimal time-varying control field $\mathbf{u}^{\star}(t)$, which is optimized for the  initial condition but converges to $\bar{\mathbf{u}}^{\star}$ over time.
}
\label{fig:tc2_dynamic}
\end{figure}

\begin{figure}[h!]
\centering
\includegraphics[width=0.5\textwidth]{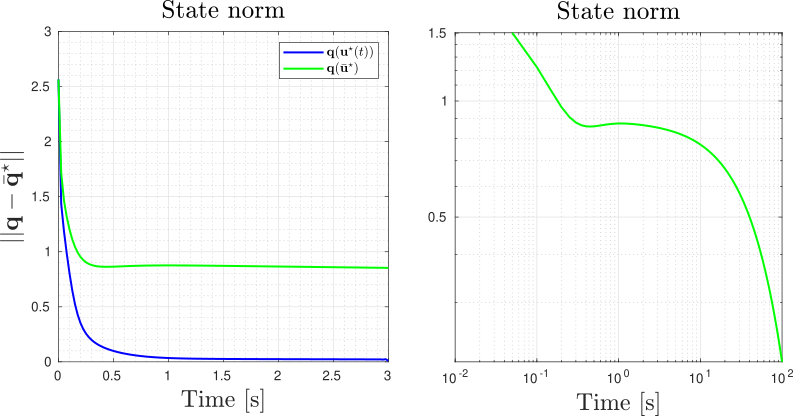}
\caption{Test case 2. Convergence in the $L^2$ norm of the density $\q$ to the optimal equilibrium $\bar{\mathbf{q}}^{\star}$ under the constant control field $\bar{\mathbf{u}}^{\star}$ (green) and the time-varying control field $\mathbf{u}^{\star}(t)$ (blue),
over ({\it left}) the first $3$ [s]; ({\it right}) $100$ [s] ($\bar{\mathbf{u}}^{\star}$ only, log-log scale). 
}
\label{fig:tc2_conv}
\end{figure}

In test case 2, 
we fabricate a much more  complicated scenario  
in which the target density has non-compact support and the domain is partially bisected by a rectangular obstacle. The optimal control field 
from the static OCP generates a unique equilibrium, as proven in Theorem \ref{stab_discrete}, but the density may be very slow to converge to this equilibrium, depending on the initial condition. Figure \ref{fig:tc2_static} plots the target density $z$, the solution $\bar{q}^{\star}$, $\bar{\mathbf{u}}^{\star}$ of the static optimization problem, and the norm of  $\bar{\mathbf{u}}^{\star}$. As in test case 1, $\bar{q}^{\star}$ is smooth, non-negative, and close to $z$, and $\bar{\mathbf{u}}^{\star}$ varies smoothly over $\Omega$ and lacks steep gradients. We solved the dynamic problem to optimize the convergence rate of the density to $\bar{q}^{\star}$ from an initial distribution that is located on the left side of the obstacle, closer to the left peak of $z$ (Fig. \ref{fig:tc2_dynamic}, top row). 
Figure \ref{fig:tc2_dynamic} presents several snapshots of the density evolution under both the constant and time-varying control fields, $\bar{\mathbf{u}}^{\star}$ and $\mathbf{u}^{\star}(t)$. The figure shows that within the first $0.6$ [s], the constant field $\bar{\mathbf{u}}^{\star}$ drives the density mostly toward the left peak of $z$, whereas the time-varying field $\mathbf{u}^{\star}(t)$ splits the density between the two peaks. Figure \ref{fig:tc2_conv} indicates that due to this splitting, $\mathbf{u}^{\star}(t)$ produces faster convergence to the optimal equilibrium density than $\bar{\mathbf{u}}^{\star}$ by two orders of magnitude; the density under  $\bar{\mathbf{u}}^{\star}$ eventually converges to the equilibrium by around $100$ [s].

\begin{figure}[h!]
\centering
\includegraphics[width=0.5\textwidth]{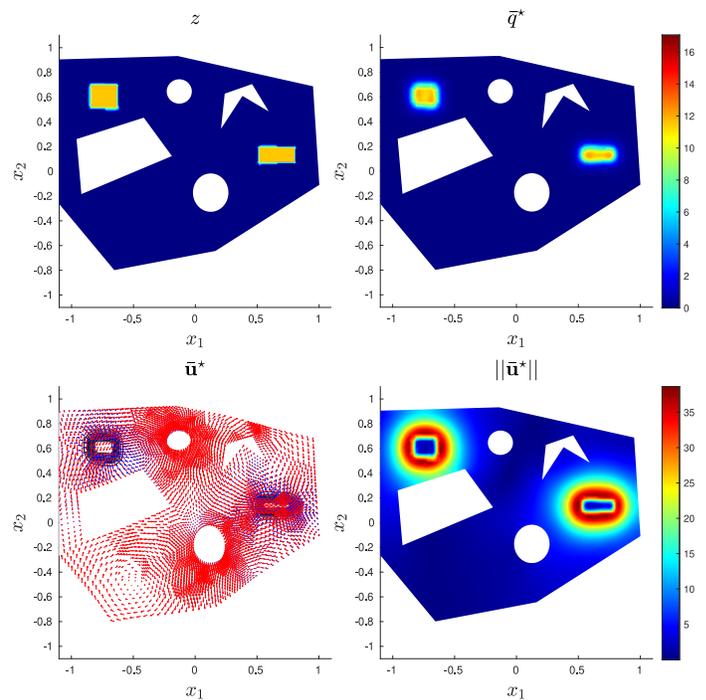}
\caption{Test case 3. Equilibrium density $\bar{q}^{\star}$ and control field $\bar{\mathbf{u}}^{\star}$ 
computed from the static optimization problem, in which $z$ is a non-smooth target function composed of the disjoint union of two characteristic functions.
The drift vector field $\mathbf{b}$ is shown in red 
in the bottom-left plot, in addition to the control vector field $\bar{\mathbf{u}}^{\star}$ in blue. The control weights 
are 
$\alpha=1$, $\beta=10^{-3}$, and $\beta_g=10^{-5}$.}
\label{fig:tc3_static}
\end{figure}

\begin{figure}[h!]
\centering
\includegraphics[width=0.4\textwidth]{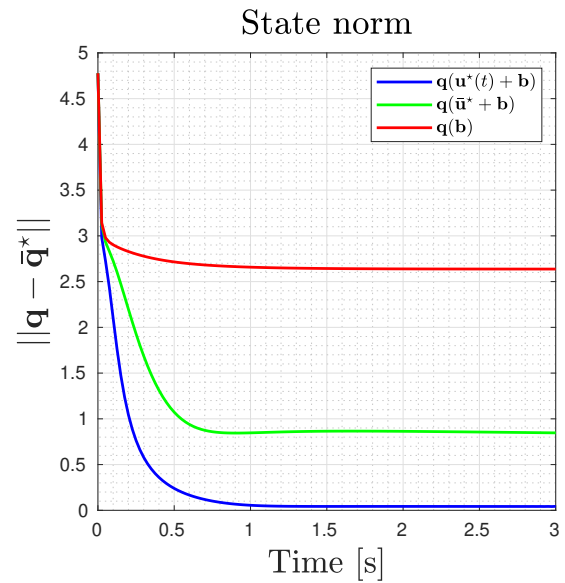}
\caption{Test case 3. $L^2$-distance between the optimal equilibrium density $\bar{\mathbf{q}}^{\star}$ and the density $\q$ under the sum of the drift velocity field $\mathbf{b}$ and (red) no control velocity field; (green) the optimal constant control field $\bar{\mathbf{u}}^{\star}$; and (blue) the optimal time-varying control field $\mathbf{u}^{\star}(t)$.}
\label{fig:tc3_conv}
\end{figure}

\begin{figure}[h!]
\centering
\includegraphics[width=0.5\textwidth]{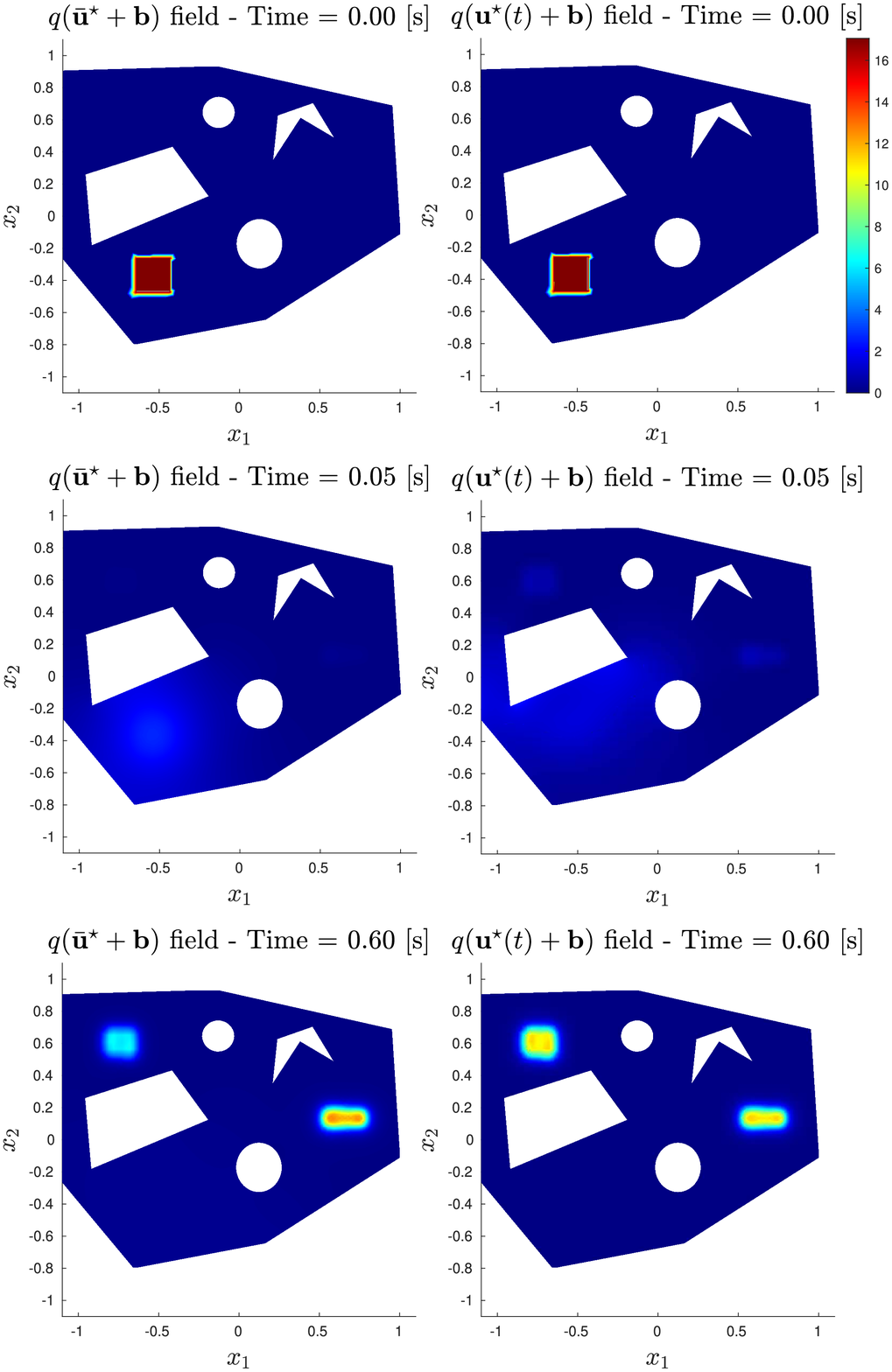}
\caption{Test case 3. Density evolution under 
the joint effect of the drift vector field $\mathbf{b}$ and ({\it left}) the optimal constant control field $\bar{\mathbf{u}}^{\star}$, or ({\it right}) the optimal time-varying control field $\mathbf{u}^{\star}(t)$, which is optimized for the initial condition but converges to $\bar{\mathbf{u}}^{\star}$ over time.}
\label{fig:tc3_3x2}
\end{figure}

In test case 3, we 
consider a scenario with an external velocity field $\mathbf{b} \in L^{\infty}(\Omega)^2$ and 
a more complex domain $\Omega$ with multiple obstacles. Introducing the drift velocity field $\mathbf{b}$ 
into the weak formulation of the state equation in the static OCP, we obtain the following problem: find $q \in \mathcal{M}_1$ such that
\begin{equation*}
    \int_{\Omega} \Big(\mu \nabla q \cdot \nabla v - (\mathbf{u}+\mathbf{b})\,q\,\nabla v \Big) \, d\Omega = 0 \quad \forall v \in H^1(\Omega).
\end{equation*}
The discretization of this state equation is: 
\begin{equation*}
   \Big( A - \mathbb{B}^{\top}\mathbf{u} - B^{\top}\Big)\mathbf{q}=\mathbf{0},
\end{equation*}
where $B$ is the transport matrix associated with 
$\mathbf{b}$, defined as $B_{ij}=\int_{\Omega} \mathbf{b} \cdot \nabla \phi_i \,\phi_j \, d\Omega$. The state equation in the dynamic OCP is modified in a similar way. Note that the additional advection term does not affect the results that we have previously derived for the OCPs.  
Figure \ref{fig:tc3_static} plots the target density $z$, the solution $\bar{q}^{\star}$, $\bar{\mathbf{u}}^{\star}$ of the static optimization problem, and the norm of  $\bar{\mathbf{u}}^{\star}$, which exhibit similar properties to the corresponding plots for the other two test cases. The drift field $\mathbf{b}$ is defined as $\mathbf{b}=[-\sin(\pi x_1)\cos(\pi x_2) ; ~\cos(\pi x_1)\sin(\pi x_2)]$. We solved the dynamic problem to obtain the time-varying control field $\mathbf{u}^{\star}(t)$ that optimizes the convergence rate of the density to $\bar{q}^{\star}$ from an initial condition 
defined as the indicator function of a square that is located at the bottom-left of the domain. Figure \ref{fig:tc3_conv} compares the convergence rate of the density to the corresponding equilibrium under $\mathbf{b}$ alone (uncontrolled), $\mathbf{b} + \bar{\mathbf{u}}^{\star}$, and $\mathbf{b} + \mathbf{u}^{\star}(t)$, and Figure \ref{fig:tc3_3x2} plots several snapshots of the density evolution under $\mathbf{b} + \bar{\mathbf{u}}^{\star}$ and $\mathbf{b} + \mathbf{u}^{\star}(t)$.
The uncontrolled density 
converges to the equilibrium induced by the drift field $\mathbf{b}$. The constant control  $\bar{\mathbf{u}}^{\star}$ 
stabilizes the optimal equilibrium density 
$\bar{\mathbf{q}}^{\star}$ in the presence of the drift field, while the time-varying 
control 
$\mathbf{u}^{\star}(t)$ speeds up the convergence to this 
equilibrium. 

%% file: Sections/conclusion.tex
In this paper, we have proposed a density control strategy 
for 
swarms of robots that follow 
single-integrator advection-diffusion dynamics. We formulated and solved an Optimal Control Problem (OCP) based on the mean-field model of the swarm to compute a space-dependent control field, defined as the robots' velocity field, that does not require inter-robot communication or density estimation algorithms for implementation. We proved that the equilibrium density of the controlled system is globally asymptotically stable, thus demonstrating that the optimal control law 
is robust to transient perturbations and independent of the initial conditions. For cases where the initial condition is approximately known, a modified dynamic OCP 
was formulated to speed up 
convergence to the optimal equilibrium density. 
Thanks to the turnpike property, the optimal solution of the dynamic OCP converges to its static counterpart, thus ensuring the stability and robustness of the control law computed by this OCP.
The analysis of the static and dynamic OCPs has been consistently carried out for both their 
infinite-dimensional formulations and their 
finite-dimensional discretizations.  
This analysis shows that several useful properties of the OCPs are inherited by the FEM discretizations.

Future work includes three main research directions. 
From a theoretical standpoint, 
we can extend our density control approach to swarms with 
explicit interactions between robots. 
This problem is addressed in \cite{Chen2021} using methods from multi-marginal optimal transport, whereas our approach would entail the formulation of a static OCP that is subject to an advection-diffusion PDE with nonlocal components in the advection 
field and is amenable to the 
stability and robustness analyses presented in this paper. 
We can also test the effectiveness of the control law in practice by implementing it on real robots, which will require incorporating 
real-world motion constraints. 
Finally, from a numerical standpoint, reduced order modeling (ROM) techniques can be used to decrease 
the computational time, and therefore bridge the gap to real-time implementation, when solutions need to be quickly recomputed as a function of certain parameters (see, e.g., \cite{sinigaglia2022fast}), such as the diffusion coefficient of the robotic swarm or the positions of obstacles in the environment.